\documentclass[10pt,a4paper]{article}
\usepackage[T1]{fontenc}
\usepackage[utf8]{inputenc}
\usepackage{lmodern}
\usepackage[margin=25mm]{geometry}
\usepackage{amsmath,amssymb,amsthm}
\usepackage{enumitem}
\usepackage{microtype}
\usepackage{hyperref}
\hypersetup{hidelinks,pdfauthor={Shuyuan Liang}}
\theoremstyle{plain}
\newtheorem{theorem}{Theorem}[section]
\newtheorem{lemma}[theorem]{Lemma}
\newtheorem{proposition}[theorem]{Proposition}
\newtheorem{corollary}[theorem]{Corollary}
\theoremstyle{definition}
\newtheorem{definition}{Definition}
\theoremstyle{remark}
\newtheorem{remark}{Remark}

\numberwithin{equation}{section}
\newcommand{\Pp}{\mathbf P}
\newcommand{\E}{\mathbf E}
\newcommand{\R}{\mathbb R}
\newcommand{\N}{\mathcal N}
\newcommand{\Phib}{\overline\Phi}
\newcommand{\eqdist}{\mathrel{\stackrel{\mathrm d}{=}}}
\DeclareMathOperator{\Var}{Var}

\newenvironment{manuscriptproof}[1][Proof]
 {\par\smallskip\noindent\textbf{#1.}\ \ignorespaces}
 {\unskip\nobreak\hfill\qed\par\smallskip}
\allowdisplaybreaks[1]

\title{Normal approximation, moderate and large deviations for the elephant random walk with stops and random step sizes}
\author{Shuyuan Liang\thanks{Corresponding author. Email:
\href{mailto:liangshuyuan0501@foxmail.com}{\nolinkurl{liangshuyuan0501@foxmail.com}}}\\[0.5ex]
\small School of Mathematics and Statistics, Northeastern University at Qinhuangdao,\\
\small Qinhuangdao, {\rm066004}, Hebei, China.}
\date{}

\begin{document}
\maketitle

\begin{abstract}
We study the elephant random walk with stops and independent positive random step sizes. Building on fixed-movement estimates, we obtain quantitative normal approximations for an empirically normalized statistic by controlling its random denominator and the independent movement clock. A permutation comparison preserves the diffusive Kolmogorov rate, while variance matching permits third-order Zolotarev estimates. We also derive self-normalized Cram\'er-type moderate deviations. For regularly varying step-size tails, we obtain uniform precise large deviations and establish the necessity of the diffusive and critical threshold orders. Further applications of the clock representation yield large deviation principles for Weibull-type and light-tailed step sizes, and a precise light-tailed expansion in the balanced case.
\end{abstract}

\begingroup\small

\par\endgroup

\section{Introduction}\label{sec:introduction}
\quad
The elephant random walk (ERW), introduced by Sch\"utz and Trimper
\cite{SchutzTrimper2004}, is a random walk whose next step depends on
its entire past. Let $X_1$ take the values $1$ and $-1$ with
$\mathbf{P}(X_1=v)=\beta\mathbf1_{\{v=1\}}+(1-\beta)\mathbf1_{\{v=-1\}},\quad v=\pm1$, where $\beta\in(0,1)$ is the bias parameter.
At time $n+1$,
choose an index $k$ uniformly from $\{1,\ldots,n\}$ and set
\[
 X_{n+1}=
 \begin{cases}
 X_k,&\text{with probability }\alpha,\\
 -X_k,&\text{with probability }1-\alpha,
 \end{cases}
 \qquad
 T_0=0,\quad T_n=\sum_{k=1}^nX_k,
\]
where $n\geq1$. Here $\alpha\in[0,1]$ is the memory parameter. The ERW is called
diffusive, critical and superdiffusive according to $\alpha<3/4$,
$\alpha=3/4$ and $\alpha>3/4$, respectively. The symmetric initial direction
gives $\E T_n=0$.

The ERW has received considerable attention. Baur and Bertoin
\cite{BaurBertoin2016} obtained functional limit theorems through
P\'olya-type urns. Coletti, Gava and Sch\"utz
\cite{Coletti2017a,Coletti2017b} established central limit theorems,
laws of large numbers and a strong invariance principle. Bercu
\cite{Bercu2018} studied the three regimes by a martingale method.
Bertoin \cite{Bertoin2022} investigated the influence of memory
on visits to the origin.

Several variants of the ERW have also been studied. Bercu and Laulin
\cite{BercuLaulin2019} investigated the multidimensional ERW,
whose recurrence and transience were studied by Qin \cite{Qin2025}.
Other extensions include ERWs with restricted memory
\cite{GutStadtmueller2021a}, smooth amnesia \cite{Laulin2022},
and polynomially decaying steps \cite{Nakano2025}.

There are two important models in these varieties of the ERW. One of them is the ERW with stops. For the ERW with stops, Cressoni, Viswanathan and da Silva
\cite{Cressoni2013} studied its moments, while Gut and Stadtm\"uller
\cite{GutStadtmueller2021b} investigated the number of zero steps.
Bercu \cite{Bercu2022} identified a Mittag--Leffler limit for the
normalized number of nonzero steps and established limit theorems
with a random normalization. The multidimensional extension was
studied in \cite{Bercu2025}. Further almost sure and Gaussian limit
results for multidimensional variants with stops or random step sizes
were obtained in \cite{GhoshDhillonKataria2026}.
In the random-step-size model of that work, zero amplitudes produce
delays with a linear asymptotic number of movements. Here the stops
are reinforced through the recalled directions, and the movement
clock has order $n^{1-r}$.
Another extension attaches independent random step sizes to the
directions. Dedecker, Fan, Hu and Merlev\`ede \cite{Dedecker2023}
obtained rates of normal approximation, and Hu and Fan
\cite{HuFan2025} proved deviation inequalities under various moment
conditions. Moment convergence rates for random step sizes were
studied by Huang, Fan, Liu and Xiang \cite{HuangFanLiuXiang2026}.
Fan and Shao \cite{FanShao2024} established Cram\'er-type
moderate deviations for bounded random step sizes, including an
empirical normalization based on squared amplitude deviations.
For related moderate and precise large deviations with random
step sizes, see also \cite{FanLiangMa2026}.
Local precise large deviations for
integer-valued step sizes were established in \cite{Liang2026}.

In this paper, we study the ERW with stops and independent positive
random step sizes. We obtain quantitative normal approximations,
self-normalized moderate deviations and large-deviation results
under finite-moment and various tail assumptions.
Our normal approximation results build on the fixed-movement bounds
of \cite{Dedecker2023}. The additional steps are the control of the
empirical denominator and the transfer through the independent
movement clock. In the diffusive regime, a permutation comparison
retains the Kolmogorov rate under empirical normalization. For
third-order Zolotarev estimates we also match the variance, whereas
the third-order Wasserstein bounds use the martingale inequality
of \cite{Guo2026}. General $p$-Wasserstein martingale bounds have
also been studied in \cite{FangKoikeSu2026}.
For heavy-tailed position deviations, we specify a uniform threshold
range and prove the necessity of its order in the diffusive and
critical regimes. The light-tailed results combine the known clock
and urn cumulants; the precise balanced-case expansion additionally
requires a local limit estimate under exponential tilting.

Section~\ref{sec:main} states the main results.
Section~\ref{sec:preliminaries} collects the preliminary estimates,
and Sections~\ref{sec:normal-proofs}--\ref{sec:large-proofs}
contain the proofs.

Throughout the paper, $c$ and $C$, possibly with indices, denote
positive constants depending only on fixed parameters and the
step-size distribution. Their values may vary from line to line.
For positive sequences, $u_n\asymp v_n$ means that their ratio is
bounded above and away from zero for all sufficiently large $n$,
and $u_n\sim v_n$ means that the ratio tends to one. We write
$u_n\gg v_n$ if $u_n/v_n\to\infty$.
Let $\Phi$ denote the standard normal distribution function and
$\Phib=1-\Phi$. 
For two random variables $U$ and $V$, $U\eqdist V$ means that $U$ and $V$ follow the same distribution.
Uniformity for $0\le x=o(h_n)$ means uniformity on
$[0,u_n]$ for every deterministic sequence $u_n=o(h_n)$.

\section{The model and main results}\label{sec:main}
\quad
The ERW with stops and random step sizes is defined as follows. Let $p,q,r\ge0$ satisfy $p+q+r=1$, with $0<r<1$. Start from the origin and
let $Y_1$ take the values $1$ and $-1$ with equal probabilities. At time
$n+1$, select an index uniformly from $\{1,\ldots,n\}$, independently of
the new randomization. If the selected increment is $Y_k$, set
\begin{equation}
 Y_{n+1}=\begin{cases}
 Y_k,&\text{with probability }p,\\
 -Y_k,&\text{with probability }q,\\
 0,&\text{with probability }r.
 \end{cases}
\end{equation}
In particular, selecting a zero increment also produces a zero. Thus $r$ is
the probability of an explicit stop in the update rule, rather than the
unconditional probability that the next displacement is zero.

Let $(Z_k)_{k\ge1}$ be a sequence of positive i.i.d.\ random variables, independent of the
entire direction and stop process, such that $\E Z_1=1$. Define
\begin{equation}\label{position-clock}
 S_n=\sum_{k=1}^n Y_kZ_k,\qquad
 \Sigma_n=\sum_{k=1}^n Y_k^2,\qquad
 b=1-r,\qquad p_r=\frac{p}{1-r}.
\end{equation}
Here $\Sigma_n$ is the number of actual movements. The memory concerns
$(Y_k)$, not the amplitudes $(Z_k)$. A positive mean different from one can
be handled by dividing the amplitudes and the position by that mean. In the
sequel, unless stated otherwise, we assume $0<p_r\le3/4$. These parameter values cover the diffusive
and critical regimes relevant to Gaussian moderate deviations.

Following the notation of the fixed-time ERW estimates, put
\begin{equation}\label{weights}
 A_m=\frac{\Gamma(m)\Gamma(2p_r)}{\Gamma(m+2p_r-1)},
 \qquad V_m=\sum_{j=1}^m A_j^2,\qquad m\ge1.
\end{equation}
In particular, $A_1=V_1=1$. Our self-normalized statistic is
\begin{equation}\label{statistic}
 W_n=\frac{A_{\Sigma_n} S_n}
 {\left[V_{\Sigma_n}
       +A_{\Sigma_n}^2\sum_{k=1}^n Y_k^2(Z_k-1)^2\right]^{1/2}}.
\end{equation}
Since $\Sigma_n\ge1$, the statistic is always defined. Its denominator
accounts for directional memory and estimates the step-size variance
$\sigma^2=\Var(Z_1)$, without requiring its value. We use self-normalization
in the empirical-variance sense of \cite[Theorem~3.2]{FanShao2024}.

\begin{remark}\label{rem:bercu-normalization}
The normalization in \eqref{statistic} differs from that used by
Bercu \cite{Bercu2022}. For unit step sizes, the latter constructs a
martingale in the original time $n$ with coefficient
\[
 \frac{\Gamma(n)\Gamma(1+p-q)}{\Gamma(n+p-q)}.
\]
Here we consider the directions at successive nonzero steps. Their
memory parameter is $p_r=p/(1-r)$, and their drift coefficient is
$2p_r-1=(p-q)/(1-r)$. Thus $A_m$ is the martingale coefficient of
the embedded ordinary ERW, evaluated at $m=\Sigma_n$ in
\eqref{statistic}. In particular, it is not obtained by replacing $n$
by $\Sigma_n$ in Bercu's coefficient.

To compare the normalizations, put $Q_m=V_m/A_m^2$. Then
\begin{equation}\label{normalizer-asymptotics}
 \begin{gathered}
 Q_m\sim
 \begin{cases}
 m/(3-4p_r),&0<p_r<3/4,\\
 m\ln m,&p_r=3/4,
 \end{cases}
 \\
 W_n=\frac{S_n}{
 [Q_{\Sigma_n}+\sum_{k=1}^nY_k^2(Z_k-1)^2]^{1/2}}.
 \end{gathered}
\end{equation}
When $Z_1=1$ almost surely, these asymptotics agree with the
normalizations in \cite[Theorems~3.3 and 3.6]{Bercu2022}; the
diffusive variance there is
$(1-r)/[3(1-r)-4p]=1/(3-4p_r)$.
For random step sizes with finite variance, the empirical term adds
$\sigma^2$ to the leading variance per movement in the diffusive
regime and is of smaller order at criticality.
All quantitative results below concern the exact statistic
\eqref{statistic}.
\end{remark}

For $0<p_r<3/4$, put
\begin{equation}
 \kappa=\min\left\{\frac12,\frac{3-4p_r}{2}\right\}.
\end{equation}
All tails below are understood separately for the two signs. They are
equal under the symmetric initial direction assumed throughout this paper.

\subsection{Rates of normal approximation}\label{subsec:normal-results}

\quad We present the normal approximation results separately for
Kolmogorov, Zolotarev and Wasserstein distances.
Let $G$ be standard normal and write
$\|U\|_s=(\E|U|^s)^{1/s}$. We first record the common weak limit
as a consequence of \cite[Theorem~3.1]{Dedecker2023}, the law of
large numbers for the empirical denominator, and the clock
representation in Lemma~\ref{lem:clock-representation}.

\begin{proposition}[A finite second moment]\label{prop:variance}
Suppose that $\E(Z_1-1)^2<\infty$ and $0<p_r\le3/4$. Then
\begin{equation}
 W_n\xrightarrow{\mathcal D}G.
\end{equation}
\end{proposition}

For the quantitative bounds, we use the moment condition
\begin{equation}\label{higher-moment}
 \E|Z_1-1|^{2+\rho}<\infty
 \quad\text{for some }0<\rho\le2.
\end{equation}

\subsubsection{Kolmogorov distance}\label{subsubsec:kolmogorov-results}

\begin{definition}[Kolmogorov distance]\label{def:kolmogorov}
For real-valued random variables $U$ and $V$, the Kolmogorov distance is
\begin{equation}\label{kolmogorov-distance}
 d_K(U,V)=\sup_{x\in\R}|\Pp(U\le x)-\Pp(V\le x)|.
\end{equation}
\end{definition}
Since the normal distribution function is continuous,
Proposition~\ref{prop:variance} also gives $d_K(W_n,G)\to0$.
For $0<p_r<3/4$, put
\begin{equation}\label{diffusive-kolmogorov-rate}
 \mathcal K_n(\rho)=n^{-b(\rho\wedge1)/2}+n^{-b(3-4p_r)}.
\end{equation}

\begin{theorem}[Kolmogorov bounds]\label{thm:kolmogorov}
The following bounds hold for $n\ge3$.
\begin{enumerate}[label=\textup{(\roman*)}]
\item If $0<p_r<3/4$ and \eqref{higher-moment} holds, then
\begin{equation}\label{low-moment-kolmogorov}
 d_K(W_n,G)\le C\mathcal K_n(\rho).
\end{equation}
For $0<p_r\le5/8$, including the balanced case $p_r=1/2$,
\begin{equation}\label{balanced-kolmogorov}
 d_K(W_n,G)\le Cn^{-b(\rho\wedge1)/2}.
\end{equation}
In particular, a finite third moment gives the rate $n^{-b/2}$
throughout this range, without an additional logarithmic factor.
\item If $p_r=3/4$ and \eqref{higher-moment} holds, then
\begin{equation}\label{kolmogorov-critical}
 d_K(W_n,G)\le\frac{C}{\ln n}.
\end{equation}
\item If $p_r=3/4$ and only $\E(Z_1-1)^2<\infty$ is assumed, then
\begin{equation}\label{kolmogorov-second-critical}
 d_K(W_n,G)\le C(\ln n)^{-1/3}.
\end{equation}
\end{enumerate}
\end{theorem}

For comparison, \cite[Theorems~3.2 and 3.3]{Dedecker2023}
already gives the corresponding Kolmogorov rate for a fixed number
of movements with deterministic normalization, including the
third-moment rate without a logarithmic loss and the $5/8$ threshold.
The diffusive assertion above retains this rate for the empirical
denominator by Lemma~\ref{lem:permutation-normal}, and then averages
over $\Sigma_n$. Thus its content concerns the normalization and
the movement clock. The next proposition records sharpness in the
balanced unit-step case by exploiting the atom at zero.

\begin{proposition}[Sharp Kolmogorov rate]\label{prop:kolmogorov-sharp}
If $Z_1=1$ almost surely and $p_r=1/2$, then
\begin{equation}\label{balanced-optimal-rate}
 d_K(W_n,G)\asymp n^{-b/2}.
\end{equation}
\end{proposition}

\subsubsection{Zolotarev distances}\label{subsubsec:zeta-results}

\begin{definition}[Zolotarev distances]\label{def:zolotarev}
Let $0<s\le3$ and let $\ell$ be the integer satisfying
$\ell<s\le\ell+1$. Let $\Lambda_s$ consist of the functions
$f\in C^\ell(\R)$ such that
$|f^{(\ell)}(x)-f^{(\ell)}(y)|\le|x-y|^{s-\ell}$ for all $x,y\in\R$.
For real-valued random variables $U$ and $V$ with finite absolute
$s$th moments, the Zolotarev distance of order $s$ is
\begin{equation}\label{probability-distances}
 \zeta_s(U,V)=\sup_{f\in\Lambda_s}|\E f(U)-\E f(V)|.
\end{equation}
For $1<s\le2$, finiteness requires equal means; for $2<s\le3$,
the second moments must also agree.
\end{definition}

Accordingly, set
\begin{equation}\label{variance-matched-statistic}
 \overline W_n=W_n/(\E W_n^2)^{1/2}.
\end{equation}
The variance is positive because $\Pp(\Sigma_n=1)>0$ and $Z_1>0$.
This correction is used only for ideal metrics of orders greater
than two. Equivalently,
$\zeta_s(W_n,(\E W_n^2)^{1/2}G)
=(\E W_n^2)^{s/2}\zeta_s(\overline W_n,G)$.

\begin{theorem}[Zolotarev bounds of orders at most two]
\label{thm:zeta-low}
Suppose that \eqref{higher-moment} holds. For $n\ge3$,
the following bounds hold.
\begin{enumerate}[label=\textup{(\roman*)}]
\item If $0<p_r<3/4$, then
\begin{equation}\label{low-moment-zeta}
 \zeta_s(W_n,G)
 \le C_s\{n^{-b(s\wedge1)\rho/(2+\rho)}
                         +n^{-b(3-4p_r)}\},
 \qquad 0<s\le2.
\end{equation}
At order two, the following estimate may be combined with
\eqref{low-moment-zeta} by taking the minimum:
\begin{equation}\label{zeta-kolmogorov-consequence}
 \zeta_2(W_n,G)\le C\mathcal K_n(\rho)\ln n.
\end{equation}
\item If $p_r=3/4$, then
\begin{equation}\label{zeta-critical}
 \zeta_s(W_n,G)\le\frac{C_s}{\ln n},\qquad 0<s\le2.
\end{equation}
\end{enumerate}
\end{theorem}

Theorem~\ref{thm:zeta-low} transfers the fixed-movement estimates
of \cite[Theorems~3.4 and 3.5]{Dedecker2023} to the empirical
normalization and random clock. Its proof combines a denominator
comparison with the clock moment bounds. The additional estimate
\eqref{zeta-kolmogorov-consequence} follows from
Theorem~\ref{thm:kolmogorov} and a weighted distribution-function
comparison.
In particular, taking $\rho=2$ gives, for $0<p_r<3/4$,
\[
 \zeta_s(W_n,G)\le
 C_s\{n^{-b(s\wedge1)/2}+n^{-b(3-4p_r)}\},
 \qquad 0<s\le2.
\]

\begin{theorem}[Zolotarev bounds under a second moment at criticality]
\label{thm:low-second-critical}
Suppose that $p_r=3/4$ and $\E(Z_1-1)^2<\infty$. Then, for $n\ge3$,
\begin{equation}\label{low-second-zeta-two}
 \zeta_2(W_n,G)\le\frac{C}{\ln n}.
\end{equation}
Consequently,
\begin{equation}\label{low-second-zeta}
 \zeta_s(W_n,G)\le C_s(\ln n)^{-s/2},\qquad 0<s<2.
\end{equation}
\end{theorem}

Here Lemma~\ref{lem:critical-comparison} supplies the comparison
with the unit-step statistic under only a second moment; the
lower-order bounds then follow by Gaussian smoothing.

\begin{proposition}[Sharp second-order Zolotarev rates]
\label{prop:normal-lower-bounds}
If $Z_1=1$ almost surely and $1/2<p_r\le3/4$, then
\begin{equation}\label{memory-lower-bounds}
 \zeta_2(W_n,G)\asymp
 \begin{cases}
 n^{-b(3-4p_r)},&1/2<p_r<3/4,\\
 (\ln n)^{-1},&p_r=3/4.
 \end{cases}
\end{equation}
\end{proposition}

These comparisons concern $W_n$, not the variance-matched
statistic $\overline W_n$.

\begin{remark}\label{rem:variance-obstruction}
For unit step sizes,
\[
 \E W_2^2=1-b\frac{(2p_r-1)^2}{1+4p_r^2}.
\]
If $p_r\ne1/2$, the variance calculation in
Section~\ref{sec:preliminaries} shows that this deficit persists for
every $n\ge2$. Thus $\zeta_3(W_n,G)=\infty$,
which explains the correction in \eqref{variance-matched-statistic}.
\end{remark}

With the variance correction in \eqref{variance-matched-statistic},
we extend the Zolotarev distance to order three. Write
\begin{equation}\label{diffusive-smooth-rate}
 R_n(\rho)=n^{-b\rho/(2+\rho)}+n^{-b(3-4p_r)},
 \qquad 0<p_r<3/4,\quad 0<\rho\le2.
\end{equation}

\begin{theorem}[Third-order Zolotarev distances]\label{thm:zeta-three}
The following assertions hold for $n\ge3$.
\begin{enumerate}[label=\textup{(\roman*)}]
\item Suppose that \eqref{higher-moment} holds with $1\le\rho\le2$.
If $0<p_r<3/4$, then
\begin{equation}\label{zeta-three-moment}
 \zeta_3(\overline W_n,G)
 \le C\left\{n^{-b\rho/(2+\rho)}+n^{-b(3-4p_r)}\right\}.
\end{equation}
In \eqref{zeta-three-moment}, finite third and fourth moments give
the respective first terms $n^{-b/3}$ and $n^{-b/2}$.
The fourth-moment assertion
includes bounded amplitudes and \eqref{exp-square}.
\item Suppose that \eqref{higher-moment} holds with $0<\rho<1$.
If $0<p_r<3/4$, then
\begin{equation}\label{zeta-three-low-moment}
 \zeta_3(\overline W_n,G)\le C R_n(\rho)\ln n.
\end{equation}
\item If $p_r=3/4$ and $\E|Z_1-1|^3<\infty$, then
\begin{equation}\label{zeta-three-critical}
 \zeta_3(\overline W_n,G)\le\frac{C}{\ln n}.
\end{equation}
\end{enumerate}
In each case the same right-hand side also bounds
$|\E W_n^2-1|$.
Moreover, throughout $0<p_r<3/4$ and $0<\rho\le2$ under
\eqref{higher-moment},
\begin{equation}\label{zeta-three-kolmogorov}
 \zeta_3(\overline W_n,G)
 \le C\mathcal K_n(\rho)(\ln n)^{3/2},
 \qquad
 |\E W_n^2-1|\le C\mathcal K_n(\rho)\ln n.
\end{equation}
The minimum of the applicable bounds may be used; no optimality
claim for these third-order bounds is intended.
\end{theorem}

The direct estimates \eqref{zeta-three-moment} and
\eqref{zeta-three-critical} use the third-order comparison in
Lemma~\ref{lem:third-smooth}, followed by variance matching.
The low-moment bound \eqref{zeta-three-low-moment} and the additional
bounds \eqref{zeta-three-kolmogorov} are consequences of the
lower-order estimates and truncation comparisons.

\begin{remark}[Intermediate Zolotarev orders]\label{rem:intermediate-orders}
For $2<s<3$, variance matching is required for $\zeta_s$ as well.
Under \eqref{higher-moment} with $1\le\rho\le2$, interpolation gives
\begin{equation}\label{intermediate-zeta}
 \zeta_s(\overline W_n,G)\le C_s
 \begin{cases}
 R_n(\rho),&0<p_r<3/4,\\
 (\ln n)^{-1},&p_r=3/4.
 \end{cases}
\end{equation}
\end{remark}

\subsubsection{Wasserstein distances}\label{subsubsec:wasserstein-results}

\begin{definition}[Wasserstein distances]\label{def:wasserstein}
For $s>0$ and real-valued random variables $U$ and $V$ with finite
absolute $s$th moments, the Wasserstein distance of order $s$ is
\begin{equation}\label{wasserstein-distance}
 \mathcal W_s(U,V)=\inf\left\{
    (\E|U'-V'|^s)^{1/(s\vee1)}:
    U'\eqdist U,\ V'\eqdist V\right\},
\end{equation}
where the infimum is over all joint laws with the given marginals.
\end{definition}

For $0<s\le1$, this transport cost equals $\zeta_s$; below we
therefore state Wasserstein results for $s\ge1$.
No variance correction is required.

\begin{theorem}[Wasserstein bounds of orders at most two]
\label{thm:wasserstein-low}
Suppose that \eqref{higher-moment} holds. For $n\ge3$,
the following bounds hold.
\begin{enumerate}[label=\textup{(\roman*)}]
\item If $0<p_r<3/4$, then
\begin{equation}\label{low-moment-wasserstein}
 \mathcal W_s(W_n,G)
 \le C_s\{n^{-b\min\{\rho/(2+\rho),\,\rho/(2s)\}}
                         +n^{-b(3-4p_r)/s}\},
 \qquad 1\le s\le2,
\end{equation}
except at $(\rho,s)=(2,2)$, where the first term is replaced by
$n^{-b/2}\sqrt{\ln n}$.
At order one, the following estimate may be combined with
\eqref{low-moment-wasserstein} by taking the minimum:
\begin{equation}\label{wasserstein-kolmogorov-consequence}
 \mathcal W_1(W_n,G)\le C\mathcal K_n(\rho)\sqrt{\ln n}.
\end{equation}
\item If $p_r=3/4$, then
\begin{equation}\label{wasserstein-critical}
 \mathcal W_s(W_n,G)\le\frac{C_s}{(\ln n)^{1/s}},
 \qquad 1\le s\le2.
\end{equation}
\end{enumerate}
\end{theorem}

Theorem~\ref{thm:wasserstein-low} combines fixed-movement transport
bounds with a coupling for the empirical denominator and
\eqref{transport-mixture}. The logarithmic factor at $(\rho,s)=(2,2)$
comes from the boundary negative moment $\E\Sigma_n^{-1}$.
The following corollary is an immediate consequence of
Theorem~\ref{thm:low-second-critical} and
\cite[Theorem~3.1]{Rio2009}.

\begin{corollary}[Wasserstein bounds under a second moment at criticality]
\label{cor:wasserstein-second-critical}
If $p_r=3/4$ and $\E(Z_1-1)^2<\infty$, then
\begin{equation}\label{wasserstein-second-critical}
 \mathcal W_s(W_n,G)\le C_s(\ln n)^{-1/2},
 \qquad 1\le s\le2,\quad n\ge3.
\end{equation}
\end{corollary}

\begin{theorem}[Third-order Wasserstein distances]\label{thm:wasserstein-three}
Suppose that $\E|Z_1-1|^3<\infty$. Then, for $n\ge3$,
\begin{equation}\label{wasserstein-three-moment}
 \mathcal W_3(W_n,G)\le C
 \begin{cases}
 n^{-b/6},&0<p_r<2/3,\\
 n^{-b/6}(\ln n)^{1/3},&p_r=2/3,\\
 n^{-b(3-4p_r)/2},&2/3<p_r<3/4,\\
 (\ln n)^{-1/2},&p_r=3/4.
 \end{cases}
\end{equation}
If, in addition, $\E|Z_1-1|^4<\infty$ and $0<p_r\le1/2$, then
\begin{equation}\label{wasserstein-three-fourth}
 \mathcal W_3(W_n,G)\le Cn^{-b/3}\sqrt{\ln n}.
\end{equation}
\end{theorem}

The lower-order estimates use the transport inequalities of
Rio \cite{Rio2009} and Bobkov \cite{Bobkov2018}.
For \eqref{wasserstein-three-moment}, we apply Guo's existing
martingale inequality \cite{Guo2026} and estimate the weight sums,
conditional variance error, and empirical denominator for this
model before averaging over the clock.
The fourth-moment bound \eqref{wasserstein-three-fourth} follows
from Theorem~\ref{thm:wasserstein-low} and a truncation comparison.
No optimality of the third-order rates is asserted.

\begin{remark}[Intermediate Wasserstein orders]
\label{rem:intermediate-wasserstein}
For $2<s<3$, the one-dimensional quantile coupling gives
\begin{equation}\label{intermediate-wasserstein}
 \mathcal W_s(W_n,G)
 \le\mathcal W_2(W_n,G)^{2(3-s)/s}
      \mathcal W_3(W_n,G)^{3(s-2)/s}.
\end{equation}
For example, a finite third moment and $0<p_r\le5/8$ give
$\mathcal W_s(W_n,G)\le C_s n^{-b/(2s)}$.
\end{remark}

\subsection{Self-normalized Cram\'er-type moderate deviations}\label{subsec:md-results}
\quad
We now turn from normal approximation in probability distances to
Cram\'er-type moderate deviations, which describe the relative
accuracy of approximating tail probabilities by the standard normal
tail \cite{Cramer1938,Petrov1975}. In the self-normalized setting,
Shao \cite{Shao1999} and Jing, Shao and Wang
\cite{JingShaoWang2003} established such estimates under finite
moment assumptions; see also \cite{DeLaPena2009}.
Related Cram\'er-type estimates for martingales were developed in
\cite{GramaHaeusler2000,FanGramaLiu2013,FanGramaLiuShao2019,FanShao2025}.

For the following estimates, assume the stronger condition
\begin{equation}\label{exp-square}
 \E\exp\{\eta(Z_1-1)^2\}<\infty
 \quad\text{for some }\eta>0.
\end{equation}
This condition is not needed for the third-moment Kolmogorov rate.
For bounded amplitudes, the fixed-movement self-normalized estimates
are given by \cite[Theorem~3.2]{FanShao2024}.
Lemma~\ref{lem:fixed-cmd} extends those estimates to
\eqref{exp-square}. Controlling small values of $\Sigma_n$ then
yields the following bounds for the stopped walk.

\begin{theorem}[Square-exponential step sizes]\label{thm:expsquare}
Suppose that \eqref{exp-square} holds. Then the following
assertions hold.
\begin{enumerate}[label=\textup{(\roman*)}]
\item If $0<p_r<3/4$, then
\begin{equation}\label{bounded-diff}
 \left|\ln\frac{\Pp(\pm W_n\ge x)}{\Phib(x)}\right|
 \le C\frac{x^3+(1+x)\ln n}{n^{b\kappa}}
\end{equation}
uniformly for $0\le x\le c\sqrt{\ln n}$, for every fixed
$0<c<\sqrt{2b(1-\kappa)}$, with a constant depending on $c$.
Moreover, on the larger range, for every fixed $0<c<\sqrt{2b}$,
\begin{equation}\label{bounded-diff-range}
 \sup_{0\le x\le c\sqrt{\ln n}}
 \left|\frac{\Pp(\pm W_n\ge x)}{\Phib(x)}-1\right|
 \longrightarrow0.
\end{equation}
\item If $p_r=3/4$, then
\begin{equation}\label{bounded-critical}
 \left|\ln\frac{\Pp(\pm W_n\ge x)}{\Phib(x)}\right|
 \le C\frac{x^3+(1+x)\ln\ln n}{\sqrt{\ln n}}
\end{equation}
uniformly for $0\le x=o(\sqrt{\ln n})$. In particular, the normal-tail
ratio tends to one uniformly for $0\le x=o((\ln n)^{1/6})$.
\end{enumerate}
\end{theorem}

In particular, Theorem~\ref{thm:expsquare} applies to bounded
positive step sizes. Condition \eqref{exp-square} also permits
unbounded step sizes and gives exponential concentration of their
empirical squared deviations. An exponential moment of $|Z_1-1|$
alone does not supply this denominator estimate.

As a direct consequence of Theorem~\ref{thm:kolmogorov} and
Mills' inequality, the normal approximation bounds also yield
relative-error estimates under finite moments.

\begin{corollary}[Finite-moment relative-error bounds]\label{cor:moment}
For all $x\ge0$ and all sufficiently large $n$, the following
assertions hold.
\begin{enumerate}[label=\textup{(\roman*)}]
\item If $0<p_r<3/4$ and \eqref{higher-moment} holds, then
\begin{equation}\label{moment-diff}
 \left|\frac{\Pp(\pm W_n\ge x)}{\Phib(x)}-1\right|
 \le C(1+x)e^{x^2/2}
 \mathcal K_n(\rho).
\end{equation}
Thus the ratio tends to one uniformly for
$0\le x\le\sqrt{2\gamma\ln n}$ whenever
\begin{equation}
 0<\gamma<b\min\left\{3-4p_r,\frac{\rho\wedge1}{2}\right\}.
\end{equation}
In particular, this holds uniformly for $x=o(\sqrt{\ln n})$.
For $0<p_r\le5/8$, the range allows every
$0<\gamma<b(\rho\wedge1)/2$.
\item If $p_r=3/4$ and \eqref{higher-moment} holds, then
\begin{equation}\label{moment-critical}
 \left|\frac{\Pp(\pm W_n\ge x)}{\Phib(x)}-1\right|
 \le C(1+x)e^{x^2/2}\frac{1}{\ln n}.
\end{equation}
The ratio therefore tends to one uniformly for
$0\le x\le\sqrt{2(1-\varepsilon)\ln\ln n}$, for every
$\varepsilon\in(0,1)$.
\item If $p_r=3/4$ and only $\E(Z_1-1)^2<\infty$ is assumed, then
\begin{equation}\label{moment-second-critical}
 \left|\frac{\Pp(\pm W_n\ge x)}{\Phib(x)}-1\right|
 \le C(1+x)e^{x^2/2}(\ln n)^{-1/3}.
\end{equation}
The ratio tends to one uniformly for
$0\le x\le\sqrt{2\gamma\ln\ln n}$ whenever $0<\gamma<1/3$.
\end{enumerate}
\end{corollary}

The factor $e^{x^2/2}$ in Corollary~\ref{cor:moment} comes from
converting an absolute normal approximation error into a relative
tail error. These estimates are not logarithmic Cram\'er expansions,
and their ranges are not claimed to be optimal. A moment of order
greater than four allows $\rho=2$.

\begin{remark}\label{rem:qualitative-md}
Under the second-moment assumption of Proposition~\ref{prop:variance},
there exists a deterministic sequence $h_n\to\infty$ such that
\begin{equation}
 \sup_{0\le x\le h_n}
 \left|\frac{\Pp(\pm W_n\ge x)}{\Phib(x)}-1\right|\longrightarrow0.
\end{equation}
Indeed, choose $h_n$ sufficiently slowly that
$d_K(W_n,G)/\Phib(h_n)\to0$, and use continuity of $\Phi$.
This consequence of the central limit theorem alone gives neither an
explicit range nor a quantitative Cram\'er expansion. At criticality,
\eqref{moment-second-critical} supplies an explicit range under the
same second-moment assumption.
\end{remark}

\subsection{Large and precise large deviations for the position $S_n$}
\label{subsec:large-results}
\quad
We finally turn to large and precise large deviations, where the
tail behavior of the random step sizes plays an important role.
For the ordinary ERW, the linear urn results of Franchini
\cite{Franchini2017,Franchini2025} provide a framework for large
deviations, while Jack and Harris \cite{JackHarris2020} studied
the underlying fluctuation mechanisms. For random step sizes,
the corresponding theory for independent sums is also relevant.
Denisov, Dieker and Shneer \cite{DenisovDiekerShneer2008}
established uniform precise asymptotics in the big-jump domain,
while Gantert, Ramanan and Rembart
\cite{GantertRamananRembart2014} obtained large deviations for
weighted sums with stretched exponential tails. In the light-tailed
setting, classical precise asymptotics were established by
Bahadur and Ranga Rao \cite{BahadurRangaRao1960} and extended
to more general sequences by Chaganty and Sethuraman
\cite{ChagantySethuraman1993}.

We now consider $S_n$ in \eqref{position-clock}. Write
$\overline F(x)=\Pp(Z_1>x)$. An LDP for $S_n/n$ with speed
$v_n\to\infty$ and good rate function $I$ means that the level sets
of $I$ are compact and, for every Borel set $B\subset\R$,
\begin{equation}\label{position-ldp}
 -\inf_{x\in B^\circ}I(x)
 \le\liminf_{n\to\infty}\frac1{v_n}\ln\Pp(S_n/n\in B)
 \le\limsup_{n\to\infty}\frac1{v_n}\ln\Pp(S_n/n\in B)
 \le-\inf_{x\in\overline B}I(x).
\end{equation}

We first establish precise large deviations for regularly varying
step sizes, uniformly over a range of thresholds. The fixed-movement
input in Lemma~\ref{lem:embedded-pld} uses the independent-sum
big-jump estimates of \cite{DenisovDiekerShneer2008}.
The random-clock argument determines the sufficient thresholds
below; Proposition~\ref{prop:heavy-threshold} subsequently addresses
their necessity in the diffusive and critical regimes.

\begin{theorem}[Regularly varying step sizes]\label{thm:position-pld}
Suppose that $0<p_r<1$ and
\begin{equation}\label{regular-tail}
 \overline F(x)=x^{-\gamma}\ell(x),\qquad \gamma>2,
\end{equation}
where $\ell$ is a slowly varying function. Define
\begin{equation}\label{position-pld-range}
 h_n^*=
 \begin{cases}
 n^{b/2}(\ln n)^{1-b/2},&0<p_r<3/4,\\
 n^{b/2}(\ln n)^{(3-b)/2},&p_r=3/4,\\
 n^{b(2p_r-1)}(\ln n)^{(2p_r-1)(1-b)+1/2},&3/4<p_r<1.
 \end{cases}
\end{equation}
Let $(h_n)_{n\ge1}$ satisfy $h_n/h_n^*\to\infty$.
Then, separately for both signs,
\begin{equation}\label{position-pld}
 \lim_{n\to\infty}\sup_{x\ge h_n}
 \left|\frac{\Pp(\pm S_n>x)}
 {\frac12\E\Sigma_n\,\overline F(x)}-1\right|=0,
 \qquad
 \E\Sigma_n=\frac{\Gamma(n+b)}{\Gamma(1+b)\Gamma(n)}.
\end{equation}
In particular, for every fixed $x>0$,
\begin{equation}\label{position-ballistic-pld}
 \Pp(\pm S_n>nx)
 \sim\frac{x^{-\gamma}}{2\Gamma(1+b)}
            n^{b-\gamma}\ell(n).
\end{equation}
\end{theorem}

The leading factor in \eqref{position-pld} is half the expected
number of movements. The following proposition shows that the
diffusive and critical threshold orders cannot be reduced.

\begin{proposition}[Necessity of the heavy-tail threshold order]
\label{prop:heavy-threshold}
Suppose that \eqref{regular-tail} holds and $0<p_r\le3/4$.
For every positive sequence $h_n=o(h_n^*)$, there exists a
sequence $x_n\ge h_n$ with $x_n=o(h_n^*)$ such that, for both signs,
\begin{equation}\label{heavy-threshold-necessity}
 \frac{\Pp(\pm S_n>x_n)}
 {\frac12\E\Sigma_n\,\overline F(x_n)}
 \longrightarrow\infty.
\end{equation}
\end{proposition}

In particular, neither logarithmic exponent in the first two lines
of \eqref{position-pld-range} can be decreased by a fixed positive
amount. The sufficiency condition $h_n/h_n^*\to\infty$ and
Proposition~\ref{prop:heavy-threshold} do not determine the sharp
constant at a fixed multiple of $h_n^*$. The superdiffusive threshold
in \eqref{position-pld-range} is only asserted to be sufficient.

For Weibull-type step-size tails, the next theorem gives logarithmic
tail asymptotics and an LDP for the position.

\begin{theorem}[Weibull-type tails]\label{thm:position-weibull}
Suppose that, for some $0<\alpha<1$ and $\beta>0$,
\begin{equation}\label{weibull-ln-tail}
 \lim_{x\to\infty}x^{-\alpha}\ln\overline F(x)=-\beta.
\end{equation}
Then, for every $x>0$,
\begin{equation}\label{position-weibull-tail}
 \lim_{n\to\infty}n^{-\alpha}\ln\Pp(\pm S_n>nx)
 =-\beta x^\alpha.
\end{equation}
Moreover, \eqref{position-ldp} holds with
\begin{equation}
 v_n=n^\alpha,\qquad I(x)=\beta|x|^\alpha.
\end{equation}
These assertions hold for every $p_r\in[0,1]$.
\end{theorem}

Condition \eqref{weibull-ln-tail} determines logarithmic
asymptotics. It does not by itself give the relative-error estimate
\eqref{position-pld} with a Weibull tail.

For the light-tailed case, put
\begin{equation}\label{clock-cumulant}
 L_b(t)=
 \begin{cases}
 -\ln[1-(1-e^{-t})^{1/b}],&t>0,\\
 0,&t\le0,
 \end{cases}
 \qquad
 J_b(s)=\sup_{t\in\R}\{st-L_b(t)\}.
\end{equation}
Let
\begin{equation}
 \lambda_{p_r}(t)=\lim_{m\to\infty}m^{-1}\ln\E e^{tT_m},
\end{equation}
where $T_m$ is the ordinary unit-step ERW with memory parameter
$0<p_r<1$. This finite, continuously differentiable function is supplied
by the linear urn results of \cite{Franchini2017,Franchini2025}.

The next theorem combines the clock cumulant from
\cite{FengHoppe1998,BercuFavaro2026} with the embedded-walk
cumulant from \cite{Franchini2017,Franchini2025}.
Conditioning on the positive and negative direction counts gives
the marked cumulant below, and the G\"artner--Ellis theorem yields
the LDP for this model.

\begin{theorem}[Light-tailed step sizes]\label{thm:position-light}
Suppose that $0<p_r<1$ and
\begin{equation}
 M(t):=\E e^{tZ_1}<\infty,\qquad t\in\R.
\end{equation}
Define
\begin{equation}\label{marked-erw-cumulant}
 \begin{aligned}
 \psi(t)&=\frac12\ln\{M(t)M(-t)\}
 +\lambda_{p_r}\!\left(\frac12\ln\frac{M(t)}{M(-t)}\right),\\
 \Lambda(t)&=L_b(\psi(t)).
 \end{aligned}
\end{equation}
Then
\begin{equation}\label{position-cumulant-limit}
 \lim_{n\to\infty}n^{-1}\ln\E e^{tS_n}=\Lambda(t),
 \qquad t\in\R,
\end{equation}
and \eqref{position-ldp} holds with
\begin{equation}\label{position-light-rate}
 v_n=n,\qquad I(x)=\sup_{t\in\R}\{tx-\Lambda(t)\}.
\end{equation}
In particular, this theorem applies to bounded step sizes and to
\eqref{exp-square}.
\end{theorem}

The composition in \eqref{marked-erw-cumulant} retains the large
deviations of the clock. Replacing $\Sigma_n$ by its expectation
does not give \eqref{position-light-rate}.

With a smooth compactly supported step-size density, we can further
obtain a precise light-tailed expansion when $p_r=1/2$.

\begin{theorem}[A precise light-tailed expansion]\label{thm:position-sharp}
Suppose that $p_r=1/2$ and $Z_1$ has a density in
$C_c^2((0,\infty))$. For $0<x<\sup\operatorname{supp}(Z_1)$,
let $t_x>0$ be the unique solution of $\Lambda'(t_x)=x$, where
\begin{equation}
 \Lambda(t)=L_b\!\left(\ln\E\cosh(tZ_1)\right).
\end{equation}
Then
\begin{equation}\label{position-sharp-tail}
 \Pp(\pm S_n\ge nx)
 \sim\frac{\exp\{-nI(x)\}}
 {b\,t_x\sqrt{2\pi n\Lambda''(t_x)}},
 \qquad I(x)=t_xx-\Lambda(t_x).
\end{equation}
The asymptotic equivalence is uniform for $x$ in every compact subset
of $(0,\sup\operatorname{supp}(Z_1))$.
\end{theorem}

The factor $b^{-1}$ in \eqref{position-sharp-tail} is inherited
from the leading term of the clock moment-generating function
in \cite{BercuFavaro2026}. The proof extends that expansion locally
to complex parameters and establishes a local limit estimate under
exponential tilting, uniformly over compact sets of tilts. This is
the model-specific step needed to obtain the precise tail
asymptotic in the framework of
\cite{BahadurRangaRao1960,ChagantySethuraman1993}.
The density assumption is a sufficient nonlattice regularity
condition for the local limit argument. Formula
\eqref{position-sharp-tail} is stated only for balanced embedded
directions; no corresponding precise expansion is asserted here
for general $p_r$.

\section{Preliminary lemmas}\label{sec:preliminaries}

\subsection{The independent random clock}
\quad
Take an ordinary ERW $(X_j)$ with memory parameter $p_r$ and
symmetric first direction, and an i.i.d. sequence
$(\widetilde Z_j)_{j\ge1}$ with $\widetilde Z_1\eqdist Z_1$.
The sequences $(X_j)$, $(\widetilde Z_j)$ and $(\Sigma_n)$
are taken to be mutually independent. Put
\begin{equation}
 \begin{gathered}
 \widetilde S_m=\sum_{j=1}^mX_j\widetilde Z_j,\qquad
 D_m^2=Q_m+\sum_{j=1}^m(\widetilde Z_j-1)^2,\\
 \widetilde W_m=\widetilde S_m/D_m.
 \end{gathered}
\end{equation}
Here $Q_m$ is defined in Remark~\ref{rem:bercu-normalization}.
For the embedded walk, write $a=2p_r-1=(p-q)/b$. Then
\begin{equation}
 Q_1=1,\qquad Q_{m+1}=(1+a/m)^2Q_m+1.
\end{equation}
Stirling's formula applied to \eqref{weights} gives
\eqref{normalizer-asymptotics}. In particular, $Q_m\ge cm$.
For later use in the superdiffusive regime, the same calculation gives
\begin{equation}\label{superdiffusive-normalizer}
 Q_m\sim q_{p_r}m^{2a},\qquad 0<q_{p_r}<\infty,
 \qquad a=2p_r-1,\quad 3/4<p_r<1.
\end{equation}
Indeed, $A_m\sim\Gamma(1+a)m^{-a}$ and
$V_m\to\sum_{j\ge1}A_j^2\in(0,\infty)$.

The positive Mittag--Leffler limit of $\Sigma_n/n^b$ was obtained
in \cite[Lemma~2.1]{Bercu2022}. The following known identity is
\cite[(4.20)]{Bercu2022}; we include its immediate moment
consequences for later reference.

\begin{lemma}[Bercu's clock moment identity]\label{lem:clock}
For every $n\ge1$ and integer $k\ge1$,
\begin{equation}\label{clock-rising-moment}
 \E\frac{\Gamma(\Sigma_n+k)}{\Gamma(\Sigma_n)}
 =\frac{k!\,\Gamma(n+bk)}{\Gamma(1+bk)\Gamma(n)}.
\end{equation}
In particular,
\begin{equation}\label{clock-positive}
 \E\Sigma_n=\frac{\Gamma(n+b)}{\Gamma(1+b)\Gamma(n)}
 \sim\frac{n^b}{\Gamma(1+b)},\qquad
 \E\Sigma_n^2\le Cn^{2b}.
\end{equation}
\end{lemma}

Identity \eqref{clock-rising-moment} is
\cite[(4.20)]{Bercu2022}; taking $k=1,2$ and using
Stirling's formula gives \eqref{clock-positive}.

We next derive the singleton probability and the negative-moment
estimates needed to control small values of the clock.

\begin{lemma}\label{lem:clock-small-values}
For every $n\ge1$,
\begin{equation}\label{clock-singleton}
 \Pp(\Sigma_n=1)=\frac{\Gamma(n-b)}{\Gamma(1-b)\Gamma(n)}
 \sim\frac{n^{-b}}{\Gamma(1-b)}.
\end{equation}
For $0<\theta<1$,
\begin{equation}\label{exact-negative}
 \E\frac{\Gamma(\Sigma_n-\theta)}{\Gamma(\Sigma_n)}
 =\frac{\Gamma(1-\theta)\Gamma(n-b\theta)}
        {\Gamma(1-b\theta)\Gamma(n)}.
\end{equation}
Moreover, for $u,v>0$,
\begin{align}
 \E\Sigma_n^{-u}
 &\le C_u
 \begin{cases}
 n^{-bu},&0<u<1,\\
 n^{-b}\ln(en),&u=1,\\
 n^{-b},&u>1,
 \end{cases}\label{clock-all-negative}\\
 \E\bigl[[\ln(e\Sigma_n)]^{-v}\bigr]
 &\le C_v[\ln(en)]^{-v}.\label{clock-logarithmic}
\end{align}
In particular, for $t>0$,
\begin{equation}\label{negative-moment}
 \Pp(\Sigma_n\le tn^b)\le C\min\{t,1\}.
\end{equation}
\end{lemma}

\begin{manuscriptproof}
The clock has transitions
\begin{equation}\label{clock-transition}
 \Pp(\Sigma_{n+1}=m+1\mid\Sigma_n=m)=bm/n,\qquad\Sigma_1=1.
\end{equation}
Hence,
\[
 \Pp(\Sigma_n=1)
 =\prod_{j=1}^{n-1}(1-b/j)
 =\frac{\Gamma(n-b)}{\Gamma(1-b)\Gamma(n)}
 \sim\frac{n^{-b}}{\Gamma(1-b)},
\]
which proves \eqref{clock-singleton}.

Identity \eqref{exact-negative} follows by taking the
Ewens--Pitman parameters $(b,-b\theta)$ in the mean-one
change-of-measure martingale of \cite{BercuFavaro2026}.

Since $\Gamma(m-u)/\Gamma(m)\asymp m^{-u}$ for $0<u<1$,
\[
 \E\Sigma_n^{-u}
 \le C_u\E\frac{\Gamma(\Sigma_n-u)}{\Gamma(\Sigma_n)}
 \le C_un^{-bu}.
\]
Thus \eqref{clock-all-negative} holds for $0<u<1$.
Write $h_{n,m}=\Pp(\Sigma_n=m)$, with zero values outside the support.
For $1\le m\le n$,
\[
 h_{n+1,m}=(1-bm/n)h_{n,m}+b(m-1)h_{n,m-1}/n,
 \qquad h_{n+1,n+1}=bh_{n,n}.
\]
Consequently, for large $n$,
\[
 \max_m h_{n+1,m}
 \le\max\{1-b/n,b\}\max_m h_{n,m}
 =(1-b/n)\max_m h_{n,m}.
\]
Iteration gives
\[
 \max_m h_{n,m}
 \le C\prod_{j=1}^{n-1}(1-b/j)
 =C\frac{\Gamma(n-b)}{\Gamma(1-b)\Gamma(n)}
 \le Cn^{-b}.
\]
Hence, for $u\ge1$,
\[
 \E\Sigma_n^{-u}
 \le Cn^{-b}\sum_{m=1}^nm^{-u}
 \le C_u
 \begin{cases}
 n^{-b}\ln(en),&u=1,\\
 n^{-b},&u>1.
 \end{cases}
\]
This completes \eqref{clock-all-negative}.

For $v>0$, splitting at $\Sigma_n=n^{b/2}$ and applying
\eqref{clock-all-negative} with $u=1/2$ gives
\[
 \E[[\ln(e\Sigma_n)]^{-v}]
 \le \Pp(\Sigma_n<n^{b/2})+C_v[\ln(en)]^{-v}
 \le Cn^{-b/4}+C_v[\ln(en)]^{-v}.
\]
This proves \eqref{clock-logarithmic}.

Finally, the point-mass bound above gives, for $t>0$,
\[
 \Pp(\Sigma_n\le tn^b)
 \le\sum_{1\le m\le\lfloor tn^b\rfloor}Cn^{-b}
 \le Ct.
\]
Combining this with the bound by one proves \eqref{negative-moment}.
\end{manuscriptproof}

The clock generating function in \cite{BercuFavaro2026}, with
$H_n(w)=\E w^{\Sigma_n}$, is
\begin{equation}\label{clock-generating-function}
 \sum_{n\ge1}H_n(w)z^{n-1}
 =\frac{(1-z)^{b-1}}{(1-z)^b-1+w^{-1}},
 \qquad w\ne0,
\end{equation}
initially for $z$ sufficiently close to zero. At $w=-1$, analytic
continuation gives
\[
 \sum_{n\ge1}\E(-1)^{\Sigma_n}z^{n-1}
 =-\frac{(1-z)^{b-1}}{2-(1-z)^b}.
\]
The denominator has no zero in $|z|\le1$. Expanding it at $z=1$
and extracting coefficients of $(1-z)^{b-1}$ therefore yields
\begin{equation}\label{clock-parity}
 \E(-1)^{\Sigma_n}\sim
 -\frac{n^{-b}}{2\Gamma(1-b)},\qquad
 \Pp(\Sigma_n\ \hbox{is even})\longrightarrow\frac12.
\end{equation}
Here the principal branch is analytic in a slit neighborhood of
the closed unit disk apart from $z=1$; the remaining singular terms
have strictly smaller coefficient order.

The next lemma gives the joint distributional representation that
transfers estimates for the embedded walk to the walk with stops.

\begin{lemma}\label{lem:clock-representation}
For every $n\ge1$, the quantities in \eqref{statistic} satisfy
the joint distributional identity
\[
 \begin{aligned}
 &\left(S_n,\sum_{k=1}^nY_k^2(Z_k-1)^2,\Sigma_n\right)\\
 &\qquad\eqdist
 \left(\widetilde S_{\Sigma_n},
              \sum_{j=1}^{\Sigma_n}(\widetilde Z_j-1)^2,\Sigma_n\right).
 \end{aligned}
\]
Consequently, we have the following decompositions 
\begin{align}
 \Pp(S_n>x)&=\sum_{m=1}^n\Pp(\Sigma_n=m)
                            \Pp(\widetilde S_m>x),\label{position-mixture}\\
 \Pp(W_n\ge x)&=\sum_{m=1}^n\Pp(\Sigma_n=m)
                            \Pp(\widetilde W_m\ge x).\label{mixture}
\end{align}
\end{lemma}

\begin{manuscriptproof}
Given movement, the recalled nonzero direction is uniform among the
previous nonzero steps and is copied with probability $p_r$.
The directions in movement order therefore have the law of an
ordinary ERW with memory parameter $p_r$, independently of the clock.
Conditional on the entire direction and stop process, the amplitudes
at successive nonzero steps have the same joint law as
$(\widetilde Z_j)_{j\ge1}$. This conditional law does not depend on
that process, proving the joint identity and the mixture formulas.
\end{manuscriptproof}

We shall use two consequences repeatedly. First,
\eqref{negative-moment} gives, for some fixed $u>0$,
\begin{equation}\label{clock-good-mass}
 \Pp(\Sigma_n\ge un^b)\ge\frac12.
\end{equation}
Secondly, Kolmogorov and Zolotarev errors can be averaged directly
over \eqref{mixture}, whereas conditional couplings give
\begin{equation}\label{transport-mixture}
 \mathcal W_s(W_n,G)^s
 \le\sum_{m=1}^n\Pp(\Sigma_n=m)
                  \mathcal W_s(\widetilde W_m,G)^s,\qquad s\ge1.
\end{equation}
To study large deviations of $S_n$, we next record the clock
cumulant and concentration bounds from
\cite{FengHoppe1998,BercuFavaro2026}, together with the
small-argument asymptotic of their rate function.

\begin{lemma}[Large-deviation estimates for the clock]
\label{lem:clock-large}
For $L_b$ and $J_b$ in \eqref{clock-cumulant},
\begin{align}
 \lim_{n\to\infty}n^{-1}\ln\E e^{t\Sigma_n}
 &=L_b(t),\qquad t\in\R,\label{clock-mgf-limit}\\
 \E e^{t\Sigma_n}&\le b^{-1}e^{nL_b(t)},\qquad t>0,
 \label{clock-mgf-upper}\\
 \Pp(\Sigma_n\ge ns)&\le b^{-1}e^{-nJ_b(s)},
 \qquad 0<s<1,\label{clock-upper-ld}\\
 J_b(s)&\sim(1-b)b^{b/(1-b)}s^{1/(1-b)},\qquad s\downarrow0.
 \label{clock-small-rate}
\end{align}
\end{lemma}

\begin{manuscriptproof}
By \eqref{clock-transition}, $\Sigma_n$ is the number of blocks
in the Ewens--Pitman model with parameters $(b,0)$.
Thus \eqref{clock-mgf-limit}--\eqref{clock-upper-ld} are the
cumulant and concentration estimates of
\cite{FengHoppe1998,BercuFavaro2026}. The moment-generating-function
bound at $n=1$ also follows from
$(1-e^{-t})^{1/b}\ge1-b^{-1}e^{-t}$.

For \eqref{clock-small-rate}, $L_b'(t)$ increases from $0$ to $1$.
As $t\downarrow0$,
$L_b(t)\sim t^{1/b}$ and $L_b'(t)\sim b^{-1}t^{(1-b)/b}$.
The maximizer $t_s$ defined by $L_b'(t_s)=s$ therefore satisfies
$t_s\sim(bs)^{b/(1-b)}$, and
\[
 J_b(s)=st_s-L_b(t_s)
 \sim(1-b)b^{b/(1-b)}s^{1/(1-b)}.
\]
\end{manuscriptproof}

\begin{lemma}[A clock-band lower bound]
\label{lem:clock-band}
There exist positive constants $c_0,C_0,C_1$, depending only on $b$,
such that, for all sufficiently large $n$ and all integers
$1\le k\le\ln n$,
\begin{equation}\label{clock-band}
 \Pp\left(c_0 n^bk^{1-b}\le\Sigma_n
                  \le C_0 n^bk^{1-b}\right)\ge e^{-C_1k}.
\end{equation}
\end{lemma}

\begin{manuscriptproof}
Put $f_k(m)=\Gamma(m+k)/\Gamma(m)$ and use the known moment
identity \eqref{clock-rising-moment}.
Uniformly for $1\le k\le2\ln n$, Stirling's formula gives
\[
 \frac{\Gamma(n+bk)}{\Gamma(n)}
 =n^{bk}\exp\{O(k^2/n)\},\qquad
 \left\{\frac{k!}{\Gamma(1+bk)}\right\}^{1/k}
 \asymp k^{1-b}.
\]
The inequalities $m^k\le f_k(m)\le(m+k)^k$ and Minkowski's
inequality consequently give
\begin{equation}\label{clock-growing-moments}
 c n^bk^{1-b}\le(\E\Sigma_n^k)^{1/k}
                  \le C n^bk^{1-b}.
\end{equation}
For the lower bound, use
$(\E\Sigma_n^k)^{1/k}\ge(\E f_k(\Sigma_n))^{1/k}-k$;
the subtracted term is uniformly negligible in this range.

Apply the Paley--Zygmund inequality to $\Sigma_n^k$. From
\eqref{clock-growing-moments} at orders $k$ and $2k$,
\[
 \Pp\left(\Sigma_n\ge
             \left[\frac12\E\Sigma_n^k\right]^{1/k}\right)
 \ge\frac{(\E\Sigma_n^k)^2}{4\E\Sigma_n^{2k}}
 \ge e^{-Ck}.
\]
This implies the same lower bound with threshold
$c_0n^bk^{1-b}$ for a sufficiently small $c_0>0$.
On the other hand, \eqref{clock-upper-ld}--\eqref{clock-small-rate}
give, for any prescribed $D>0$ and a sufficiently large $C_0$,
\[
 \Pp(\Sigma_n>C_0n^bk^{1-b})\le e^{-Dk}
\]
uniformly for $1\le k\le\ln n$ and large $n$.
Taking $D$ larger than the constant in the lower bound, subtracting
this upper tail, and enlarging $C_1$ proves \eqref{clock-band}.
\end{manuscriptproof}

\subsection{Moment and normalization estimates}

Put
\[
 \varepsilon_1=X_1,\qquad
 \varepsilon_j=X_j-a\frac{T_{j-1}}{j-1},\quad j\ge2.
\]
These are martingale differences relative to
$\mathcal F_j=\sigma(X_1,\ldots,X_j)$, with $|\varepsilon_k|\leq2,\quad 1\leq k\leq n$. Since $A_mT_m=\sum_{j=1}^m A_j\varepsilon_j$, the
conditional Hoeffding inequality, obtained by applying
\cite[(4.16)]{Hoeffding1963} conditionally and iterating, gives
\begin{equation}\label{unit-subgaussian}
 \E\exp\{tT_m/\sqrt{Q_m}\}\le e^{t^2/2},
 \qquad \|T_m\|_s\le C_s\sqrt{Q_m},\quad s\ge1.
\end{equation}
These inequalities hold throughout $0<p_r<1$; no central limit
theorem is needed for \eqref{unit-subgaussian}.

In the diffusive regime, the following bound for weighted sums
of the directions also controls the conditional moments of
$\widetilde W_m$.

\begin{lemma}\label{lem:weighted-directions}
If $0<p_r<3/4$, then for deterministic real $c_1,\ldots,c_m$,
\begin{equation}
 \E\exp\left\{t\sum_{j=1}^m c_jX_j\right\}
 \le\exp\left\{C_at^2\sum_{j=1}^m c_j^2\right\}.
\end{equation}
Consequently, without a moment condition on the amplitudes, we have
\begin{equation}\label{conditional-self-moments}
 \E\left(e^{t\widetilde W_m}\mid
                 \widetilde Z_1,\ldots,\widetilde Z_m\right)
 \le e^{Ct^2},
 \qquad
 \E\left(|\widetilde W_m|^s\mid
                 \widetilde Z_1,\ldots,\widetilde Z_m\right)
 \le C_s,\quad s\ge1.
\end{equation}
\end{lemma}

\begin{manuscriptproof}
The martingale representation yields
\[
 \begin{gathered}
 \sum_{j=1}^m c_jX_j=\sum_{j=1}^m d_j\varepsilon_j,\\
 d_j=c_j+aA_j\sum_{k=j+1}^m
       \frac{c_k}{(k-1)A_{k-1}},\quad 1\le j<m,
 \qquad d_m=c_m.
 \end{gathered}
\]
By \eqref{weights}, $A_j\asymp j^{-a}$, and hence
\[
 \left|\frac{aA_j}{(k-1)A_{k-1}}\right|
 \le C_aj^{-a}k^{a-1},\qquad 1\le j<k\le m.
\]
For Schur's test, take $w_j=j^{-1/2}$. Since $a<1/2$, the
weighted row and column sums satisfy
\begin{align*}
 j^{-a}\sum_{k=j+1}^m k^{a-1}w_k
 &=j^{-a}\sum_{k=j+1}^m k^{a-3/2}
 \le j^{-a}\int_j^\infty x^{a-3/2}\,dx
 =\frac{w_j}{1/2-a},\\
 k^{a-1}\sum_{j=1}^{k-1}j^{-a}w_j
 &=k^{a-1}\sum_{j=1}^{k-1}j^{-a-1/2}
 \le C_ak^{a-1}k^{1/2-a}=C_aw_k.
\end{align*}
The corresponding weighted Cauchy--Schwarz estimate gives
\begin{align*}
 \sum_{j=1}^m(d_j-c_j)^2
 &\le C_a\sum_{j=1}^{m-1}
       \left(\sum_{k=j+1}^m j^{-a}k^{a-1}|c_k|\right)^2\\
 &\le C_a\sum_{j=1}^{m-1}
       \left(\sum_{k=j+1}^m j^{-a}k^{a-1}w_k\right)
       \left(\sum_{k=j+1}^m
                    j^{-a}k^{a-1}\frac{c_k^2}{w_k}\right)\\
 &\le C_a\sum_{j=1}^{m-1}w_j
          \sum_{k=j+1}^m j^{-a}k^{a-1}\frac{c_k^2}{w_k}\\
 &=C_a\sum_{k=2}^m\frac{c_k^2}{w_k}\,
          k^{a-1}\sum_{j=1}^{k-1}j^{-a}w_j
 \le C_a\sum_{k=1}^m c_k^2.
\end{align*}
Consequently,
\[
 \sum_{j=1}^m d_j^2
 \le2\sum_{j=1}^m c_j^2+2\sum_{j=1}^m(d_j-c_j)^2
 \le C_a\sum_{j=1}^m c_j^2.
\]
Set $\mu_1=0$ and
$\mu_j=aT_{j-1}/(j-1)$ for $j\ge2$, and take $\mathcal F_0$ to be
trivial. Conditional on $\mathcal F_{j-1}$, the difference
$\varepsilon_j$ is centered and lies in $[-1-\mu_j,1-\mu_j]$,
an interval of length $2$. Applying Hoeffding's lemma
\cite[(4.16)]{Hoeffding1963} to this conditional distribution gives
\[
 \E\left(e^{td_j\varepsilon_j}\mid\mathcal F_{j-1}\right)
 \le e^{t^2d_j^2/2},\qquad t\in\R.
\]
Since $d_j$ is deterministic, conditional iteration and
$\sum_jd_j^2\le C_a\sum_jc_j^2$ give
\[
 \E\exp\left\{t\sum_{j=1}^m c_jX_j\right\}
 =\E\exp\left\{t\sum_{j=1}^m d_j\varepsilon_j\right\}
 \le\exp\left\{\frac{t^2}{2}\sum_{j=1}^m d_j^2\right\}
 \le\exp\left\{C_at^2\sum_{j=1}^m c_j^2\right\}.
\]
By independence of the directions and amplitudes, we may condition on
$(\widetilde Z_1,\ldots,\widetilde Z_m)$ and apply the preceding bound
with $c_j=\widetilde Z_j/D_m$. Since $Q_m\ge cm$,
\[
 \sum c_j^2
 \le\frac{2m+2\sum(\widetilde Z_j-1)^2}
              {Q_m+\sum(\widetilde Z_j-1)^2}\le C.
\]
This proves the conditional exponential moment bound in
\eqref{conditional-self-moments}. The conditional moment bounds follow
by integrating the corresponding two-sided exponential tail bound.
\end{manuscriptproof}

The same directional martingale gives, for $m\ge2$,
\begin{equation}\label{unit-variance-deficit}
 1-\frac{\E T_m^2}{Q_m}
 =\frac{a^2}{V_m}\sum_{j=2}^m
      \frac{A_j^2\E T_{j-1}^2}{(j-1)^2}.
\end{equation}
If $a\ne0$, the right-hand side is positive. For unit amplitudes,
averaging over the clock therefore gives $\E W_n^2<1$ for every
$n\ge2$, establishing the obstruction in
Remark~\ref{rem:variance-obstruction}.

Whenever $\sigma^2<\infty$, put
\begin{equation}
 B_m^2=Q_m+m\sigma^2,\qquad U_m=\widetilde S_m/B_m,\qquad
 \Delta_m=\sum_{j=1}^m[(\widetilde Z_j-1)^2-\sigma^2].
\end{equation}
Thus $D_m^2=B_m^2+\Delta_m$. The decomposition
\begin{equation}
 \widetilde S_m=T_m+\sum_{j=1}^mX_j(\widetilde Z_j-1)
\end{equation}
will be used throughout. Conditional on the directions, the second
term is a sum of independent centered variables. Symmetrization and
the Rademacher moment inequality give
\begin{equation}\label{amplitude-noise-moments}
 \E(\widetilde S_m-T_m)^2=m\sigma^2,\qquad
 \|\widetilde S_m-T_m\|_s\le C_s\sqrt m
 \quad\text{if }s\ge2,\ \E|Z_1-1|^s<\infty.
\end{equation}
In particular, $\|U_m\|_s\le C_s$ under the same assumption.

For each fixed $m$, reveal the directions and then the centered
amplitudes. This gives a martingale with terminal value $U_m$ and
successive differences
\begin{equation}\label{two-stage}
 \frac{A_j\varepsilon_j}{A_mB_m},\quad 1\le j\le m,
 \qquad
 \frac{X_j(\widetilde Z_j-1)}{B_m},\quad 1\le j\le m.
\end{equation}
Its conditional variance sum is
\begin{equation}\label{two-stage-bracket}
 1-\frac{a^2}{A_m^2B_m^2}
       \sum_{j=2}^m\frac{A_j^2T_{j-1}^2}{(j-1)^2}.
\end{equation}
The construction is for a fixed terminal time; the process
$(A_m\widetilde S_m)$ itself need not be a martingale.

Under \eqref{higher-moment}, the von Bahr--Esseen inequality
\cite{VonBahrEsseen1965} yields
\begin{equation}\label{square-moment}
 \E|\Delta_m|^{1+\rho/2}\le Cm.
\end{equation}
Since $B_m^2\ge cm$, Markov's inequality gives
\begin{equation}
 \Pp\left(\left|D_m^2/B_m^2-1\right|>\varepsilon\right)
 \le Cm^{-\rho/2}\varepsilon^{-(1+\rho/2)},
 \qquad 0<\varepsilon<\tfrac12.
\end{equation}
Under \eqref{exp-square}, the logarithm of the moment generating
function of $(Z_1-1)^2-\sigma^2$ is bounded by $C\lambda^2$ for
$|\lambda|\le c$. Exponential Markov inequality therefore gives
\begin{equation}\label{squares-concentration}
 \Pp(|\Delta_m|>mt)\le2\exp\{-cm\min(t^2,t)\},\qquad t>0.
\end{equation}
We also have the deterministic bounds
\begin{equation}
 D_m/B_m\ge c,\qquad
 |B_m/D_m-1|\le C\min\{|\Delta_m|/B_m^2,1\}.
\end{equation}

The following smooth comparison measures the contribution of the
random step sizes and will be used in the critical regime.

\begin{lemma}\label{lem:critical-comparison}
Let $s\in\{2,3\}$ and assume that $\E|Z_1-1|^s<\infty$.
Uniformly over $f\in\Lambda_s$ with $f^{(j)}(0)=0$ for
$0\le j<s$,
\begin{equation}
 |\E f(\widetilde W_m)-\E f(T_m/\sqrt{Q_m})|
 \le C_s\left\{\frac{m}{Q_m}
                    +\left(\frac{m}{Q_m}\right)^{s/2}\right\}.
\end{equation}
\end{lemma}

\begin{manuscriptproof}
For the stated test functions,
\[
 |f(x+y)-f(x)-f'(x)y|
 \le C_s\{|x|^{s-2}y^2+|y|^s\},\qquad
 \E(\widetilde S_m-T_m\mid X_1,\ldots,X_m)=0.
\]
Hence \eqref{unit-subgaussian} and
\eqref{amplitude-noise-moments} imply
\begin{align*}
 &\left|\E f(\widetilde S_m/\sqrt{Q_m})
                -\E f(T_m/\sqrt{Q_m})\right|\\
 &\quad\le C_s Q_m^{-s/2}
  \{m\sigma^2\E|T_m|^{s-2}+\E|\widetilde S_m-T_m|^s\}
 \le C_s\{m/Q_m+(m/Q_m)^{s/2}\}.
\end{align*}
Moreover,
\begin{align*}
 &|f(\widetilde S_m/\sqrt{Q_m})-f(\widetilde S_m/D_m)|\\
 &\quad\le C_s Q_m^{-s/2}
   (|T_m|^s+|\widetilde S_m-T_m|^s)
        \min\{(D_m^2-Q_m)/Q_m,1\}.
\end{align*}
Using $T_m\perp(D_m^2-Q_m)$ and
$\E(D_m^2-Q_m)=m\sigma^2$ gives the same bound in expectation.
\end{manuscriptproof}

\subsection{Gaussian tail estimates}

For bounded amplitudes, the estimates in the following lemma,
including the empirical denominator, are those of Fan and Shao
\cite[Theorem~3.2]{FanShao2024}. We give the extension to the
square-exponential moment condition by verifying the martingale
conditions and controlling the empirical squared deviations.

\begin{lemma}\label{lem:fixed-cmd}
Suppose \eqref{exp-square} holds. If $0<p_r<3/4$, then
\begin{equation}\label{fixed-diff}
 \left|\ln\frac{\Pp(\pm\widetilde W_m\ge x)}{\Phib(x)}\right|
 \le C\{x^3+(1+x)\ln m\}m^{-\kappa}
\end{equation}
uniformly for $0\le x=o(m^\kappa)$. If $p_r=3/4$, then
\begin{equation}\label{fixed-critical}
 \left|\ln\frac{\Pp(\pm\widetilde W_m\ge x)}{\Phib(x)}\right|
 \le C\frac{x^3+(1+x)\ln\ln m}{\sqrt{\ln m}}
\end{equation}
uniformly for $0\le x=o(\sqrt{\ln m})$.
\end{lemma}

\begin{manuscriptproof}
Condition \eqref{exp-square} implies, for some $H\ge1$,
\[
 \E|Z_1-1|^\ell\le\frac12\ell!H^{\ell-2}\sigma^2,
 \qquad \ell\ge2.
\]
Let $(\xi_j,\mathcal H_j)_{1\le j\le2m}$ denote the differences
and filtration in \eqref{two-stage}, and put
\[
 \delta_m=
 \begin{cases}
 m^{-\kappa},&p_r<3/4,\\
 [\ln(em)]^{-1/2},&p_r=3/4.
 \end{cases}
\]
Since $|\varepsilon_j|\le2$ and
\[
 \max_{1\le j\le m}\frac{A_j}{A_mB_m}+\frac{H}{B_m}
 \le C\delta_m,
\]
we have, for $\ell\ge2$ and $1\le j\le2m$,
\[
 \E(|\xi_j|^\ell\mid\mathcal H_{j-1})
 \le\frac12\ell!(C\delta_m)^{\ell-2}
                  \E(\xi_j^2\mid\mathcal H_{j-1}).
\]
Moreover, \eqref{unit-subgaussian} and the weight asymptotics give
\[
 \sup_{j\ge1}\E\exp\{cT_j^2/Q_j\}\le C,\qquad
 \frac{a^2}{A_m^2B_m^2}
 \sum_{j=2}^m\frac{A_j^2Q_{j-1}}{(j-1)^2}
 \le C\delta_m^2.
\]
By \eqref{two-stage-bracket} and Jensen's inequality,
\[
 \E\exp\left\{\frac{c}{\delta_m^2}
 \left(1-\sum_{j=1}^{2m}
       \E(\xi_j^2\mid\mathcal H_{j-1})\right)\right\}\le C.
\]
The variance deficit is nonnegative; hence, for $u>0$,
\[
 \Pp\left(\left|\sum_{j=1}^{2m}
       \E(\xi_j^2\mid\mathcal H_{j-1})-1\right|>u\right)
 \le Ce^{-cu/\delta_m^2}.
\]
Thus \cite[Theorem~2.2 and Remark~2.1]{FanShao2024} yields
\eqref{fixed-diff}--\eqref{fixed-critical} with $U_m$ in place
of $\widetilde W_m$. The argument also includes $\sigma^2=0$.

Choose a sufficiently large constant $D$ and set
\[
 \varepsilon_m=\frac D{\sqrt m}
 \begin{cases}
 x+\sqrt{\ln m},&p_r<3/4,\\
 x+\sqrt{\ln\ln m},&p_r=3/4.
 \end{cases}
\]
On the respective ranges, $\varepsilon_m=o(1)$ and
\begin{equation}
 |D_m^2/B_m^2-1|\le\varepsilon_m
\end{equation}
except on an event of probability at most $q_m$, where
\[
 q_m\le C e^{-Kx^2}
 \begin{cases}m^{-K},&p_r<3/4,\\(\ln m)^{-K},&p_r=3/4.
 \end{cases}
\]
Here $K$ can be arbitrarily large, by \eqref{squares-concentration}.
The tail sandwich and Mills' inequality give
\begin{gather*}
 \Pp(U_m\ge x\sqrt{1+\varepsilon_m})-q_m
 \le\Pp(\widetilde W_m\ge x)
 \le\Pp(U_m\ge x\sqrt{1-\varepsilon_m})+q_m,\\
 \frac{q_m}{\Phib(x)}\le C(1+x)e^{x^2/2}q_m,
\end{gather*}
and
\begin{equation}
 \left|\ln\frac{\Phib(x\sqrt{1+v})}{\Phib(x)}\right|
 \le C(1+x^2)|v|,\qquad |v|\le\tfrac12,\quad x\ge0.
\end{equation}
Furthermore,
\[
 (1+x^2)\varepsilon_m\le C
 \begin{cases}
 \{x^3+(1+x)\ln m\}m^{-\kappa},&p_r<3/4,\\
 \{x^3+(1+x)\ln\ln m\}/\sqrt{\ln m},&p_r=3/4.
 \end{cases}
\]
The logarithmic errors for $U_m$ are $o(x^2)+o(1)$ on the
respective ranges. Thus $q_m$ is negligible relative to the
lower tail bound and, when the stated error tends to zero,
relative to that error times $\Phib(x)$. The sandwich proves
both estimates. Reflection gives the other tail.
\end{manuscriptproof}

\section{Proofs of the normal approximation results}\label{sec:normal-proofs}

\subsection{The central limit theorem and lower-order bounds}

\begin{manuscriptproof}[Proof of Proposition~\ref{prop:variance}]
The fixed-movement central limit theorem
\cite[Theorem~3.1]{Dedecker2023}, with symmetric initial direction,
gives $U_m\Rightarrow G$, including when $\sigma^2=0$.
The law of large numbers yields $D_m^2/B_m^2\to1$ in probability,
so $\widetilde W_m\Rightarrow G$ by Slutsky's theorem.
Since $\Sigma_n\to\infty$ in probability, \eqref{mixture} gives
$W_n\Rightarrow G$.
\end{manuscriptproof}

We first recall the fixed-movement inputs, with the symmetric initial
direction used in this paper. Theorems~3.2 and 3.4 of
\cite{Dedecker2023} give
\begin{align}
 d_K(T_m/\sqrt{Q_m},G)
 &\le C(m^{-1/2}+V_m^{-1}),\label{unit-low-kolmogorov}\\
 \zeta_s(T_m/\sqrt{Q_m},G)
 &\le C_s(m^{-s/2}+V_m^{-1}),\qquad 0<s\le2.
 \label{unit-low-zeta}
\end{align}
For the deterministic normalization $U_m$ defined above,
Theorems~3.3 and 3.5 of \cite{Dedecker2023}, with the moment
parameter reduced to $\rho\wedge1$ when necessary, yield
\begin{align}
 d_K(U_m,G)&\le C\{m^{-(\rho\wedge1)/2}+V_m^{-1}\},
 \\
 \zeta_s(U_m,G)&\le
 C_s\{m^{-(s\wedge\rho\wedge1)/2}+V_m^{-1}\},
 \qquad 0<s\le2.\label{deterministic-low-zeta}
\end{align}
For $1\le s\le2$ we shall also use
\begin{equation}\label{deterministic-low-wasserstein}
 \mathcal W_s(U_m,G)
 \le C_s\{m^{-(\rho\wedge s)/(2s)}+V_m^{-1/s}\}.
\end{equation}
For $0<\rho\le1$, this follows from
\eqref{deterministic-low-zeta} and
\cite[Theorem~3.1]{Rio2009}, with $\mathcal W_1=\zeta_1$.
For $1<\rho\le2$ and $\sigma^2>0$,
\cite[Corollary~1.2]{Bobkov2018}, applied conditionally to the
independent variables $X_j(\widetilde Z_j-1)$, gives
\[
 \sup_{(X_j)}\mathcal W_s\!\left(
 \mathcal L\left(\frac{\sum_jX_j(\widetilde Z_j-1)}{\sigma\sqrt m}
                      \,\middle|\,X_1,\ldots,X_m\right),\N(0,1)\right)
 \le C_s m^{-(\rho\wedge s)/(2s)}.
\]
Combining this with \eqref{unit-low-zeta} by conditional Gaussian
couplings proves \eqref{deterministic-low-wasserstein}.
For $\sigma^2=0$, only the unit-step estimate is needed.

\begin{lemma}[Permutation comparison]\label{lem:permutation-normal}
Suppose that $0<p_r<3/4$ and \eqref{higher-moment} holds.
With $\delta=\rho\wedge1$, for every $m\ge1$,
\begin{equation}\label{fixed-permutation-kolmogorov}
 d_K(\widetilde W_m,G)
 \le C\{m^{-\delta/2}+V_m^{-1}\}.
\end{equation}
\end{lemma}

\begin{manuscriptproof}
If $\sigma^2=0$, then $Z_1=1$ almost surely and
\eqref{unit-low-kolmogorov} proves the assertion. Suppose
$\sigma^2>0$, and take $m\ge4$, absorbing smaller $m$ into $C$.
Write $z_j=\widetilde Z_j$,
\[
 \bar z=m^{-1}\sum_{j=1}^m z_j,\qquad
 s_m^2=m^{-1}\sum_{j=1}^m(z_j-\bar z)^2,\qquad
 D^2=Q_m+\sum_{j=1}^m(z_j-1)^2.
\]
Let $\pi$ be a uniform random permutation, independent of
$(T_m,z_1,\ldots,z_m)$, and put
\[
 B_{m,k}=2\sum_{j=1}^kz_{\pi(j)}-\sum_{j=1}^mz_j,\qquad0\le k\le m.
\]
Independence and identical distribution of the amplitudes imply
\[
 \widetilde W_m\eqdist B_{m,(m+T_m)/2}/D.
\]
Indeed, conditional on any direction vector with $k$ positive
coordinates, the joint law of the numerator and the symmetric
denominator depends only on $k$. This identity does not require the
direction sequence to be exchangeable.

Let $H_m\sim\N(0,Q_m)$ be independent of the amplitudes and
permutation, and let $L_m$ be its nearest point in
$\{-m,-m+2,\ldots,m\}$, with clipping at the endpoints.
Since $Q_m\asymp m$, \eqref{unit-low-kolmogorov} and Gaussian
anti-concentration give
\[
 d_K(T_m,L_m)\le C(m^{-1/2}+V_m^{-1}),\qquad
 \Pp(|L_m|>m/2)\le Ce^{-cm},\qquad \E L_m^2\le Cm.
\]
For each fixed $(z,\pi)$, $B_{m,k}/D$ increases with $k$ because
$z_j>0$. Thus replacing $T_m$ by $L_m$ in the representation costs
at most $C(m^{-1/2}+V_m^{-1})$ in Kolmogorov distance.

Conditional on $z$ and $L_m=t=2k-m$, the permutation mean and
variance are
\begin{equation}\label{permutation-variance}
 \E_\pi B_{m,k}=t\bar z,\qquad
 v_{m,t}:=\mathop{\rm Var}\nolimits_\pi(B_{m,k})
       =\frac{m^2-t^2}{m-1}s_m^2.
\end{equation}
For $|t|\le m/2$ and $s_m>0$, use the deterministic array
\[
 u_i=\mathbf1_{\{i\le k\}}-\frac{k}{m},\qquad
 c_{ij}=\frac{2u_i(z_j-\bar z)}{\sqrt{v_{m,t}}}.
\]
Its row and column sums vanish and
$(m-1)^{-1}\sum_{i,j}c_{ij}^2=1$. Applying
\cite[Theorem~1.1 and (1.5)]{ChenFang2015} therefore gives
\begin{equation}\label{permutation-be}
 d_K\left(\frac{B_{m,k}-t\bar z}{\sqrt{v_{m,t}}},G\right)
 \le\min\left\{1,\frac{451}{m}\sum_{i,j}|c_{ij}|^3\right\}
 \le\min\left\{1,\frac{C\sum_j|z_j-\bar z|^3}
                         {m^{3/2}s_m^3}\right\}.
\end{equation}
The distance in this display concerns only the permutation
randomness, with $z$ and $t$ fixed.

There exists $v_0>0$ such that
\begin{equation}\label{sample-variance-lower}
 \Pp(s_m^2<v_0)\le Ce^{-cm}.
\end{equation}
To see this, choose two bounded intervals at positive distance,
each having positive $Z_1$-probability. Except on an exponentially
small binomial event, each contains a fixed positive fraction of
the sample. The identity
$s_m^2=(2m^2)^{-1}\sum_{i,j}(z_i-z_j)^2$ then gives
\eqref{sample-variance-lower}.
Also,
\[
 \Pp(\max_jz_j>\sqrt m)\le Cm^{-\delta/2},\qquad
 \E\sum_j|z_j-\bar z|^{2+\delta}\le Cm.
\]
On $\{\max_jz_j\le\sqrt m\}$, positivity gives
$|z_j-\bar z|\le\sqrt m$, and hence
\[
 \sum_j|z_j-\bar z|^3
 \le m^{(1-\delta)/2}\sum_j|z_j-\bar z|^{2+\delta}.
\]
Consequently, the expectation of the right-hand side of
\eqref{permutation-be}, restricted to $s_m^2\ge v_0$, is at most
$Cm^{-\delta/2}$. The omitted events contribute the same order.
With $G_2$ an additional independent standard normal variable,
we have proved
\[
 d_K\left(\widetilde W_m,
       \frac{L_m\bar z+\sqrt{v_{m,L_m}}G_2}{D}\right)
 \le C\{m^{-\delta/2}+V_m^{-1}\}.
\]

On $|L_m|\le m/2$ and $s_m^2\ge v_0$, the relative variance
difference between $v_{m,L_m}$ and $ms_m^2$ is at most
$C(m^{-1}+L_m^2/m^2)$. Gaussian scaling and $\E L_m^2\le Cm$
show that replacing $\sqrt{v_{m,L_m}}G_2$ by $\sqrt m s_mG_2$
costs $O(m^{-1})$, after averaging. Next replace $L_m$ by $H_m$.
Outside an exponentially small clipping event, $|H_m-L_m|\le1$.
The normal density bound, conditional on $H_m$ and $z$, bounds
the resulting distance by $C\bar z/(\sqrt m s_m)$ on
$s_m^2\ge v_0$. Since $\E\bar z=1$, its average is
$O(m^{-1/2})$.

Conditional on $z$, the resulting variable is centered normal with
variance
\[
 r_m=\frac{Q_m\bar z^2+ms_m^2}{D^2},\qquad
 D^2=Q_m+ms_m^2+m(\bar z-1)^2.
\]
The cancellation
\[
 Q_m\bar z^2+ms_m^2-D^2
 =2Q_m(\bar z-1)+(Q_m-m)(\bar z-1)^2
\]
and $D^2\ge Q_m\asymp m$ imply
\[
 \E|r_m-1|
 \le C\E\{|\bar z-1|+(\bar z-1)^2\}\le Cm^{-1/2}.
\]
Finally, $d_K(\sqrt vG,G)\le C\min\{1,|v-1|\}$ for $v>0$.
Combining the comparisons proves
\eqref{fixed-permutation-kolmogorov}.
\end{manuscriptproof}

\begin{manuscriptproof}[Proof of Theorem~\ref{thm:kolmogorov}]
In the diffusive regime, Lemma~\ref{lem:permutation-normal}
provides the fixed-movement estimate. Average it using
\eqref{mixture} and \eqref{clock-all-negative}.
Writing $\delta=\rho\wedge1$ and $\alpha=3-4p_r$, we have
$\delta/2<1$. If $\alpha<1$, then
$\E V_{\Sigma_n}^{-1}\le Cn^{-b\alpha}$.
If $\alpha\ge1$, this contribution is at most $Cn^{-b}\ln(en)$
and is absorbed by $n^{-b\delta/2}$.
For $p_r\le5/8$, $\alpha\ge1/2\ge\delta/2$, proving
\eqref{balanced-kolmogorov}.

At criticality under \eqref{higher-moment}, put
$g_m=|B_m/D_m-1|$. Then
\[
 g_m\le C,\qquad \E g_m^{1+\rho/2}\le Cm^{-\rho/2},
 \qquad \|U_m\|_{2+\rho}\le C.
\]
H\"older's inequality yields the comparison
\begin{equation}\label{critical-normalization-l2}
 \|\widetilde W_m-U_m\|_2^2
 \le C\{\E g_m^{2(2+\rho)/\rho}\}^{\rho/(2+\rho)}
 \le Cm^{-\rho^2/(2(2+\rho))}.
\end{equation}
The elementary coupling inequality gives
\[
 d_K(\widetilde W_m,G)\le d_K(U_m,G)+C\varepsilon
           +\varepsilon^{-2}\|\widetilde W_m-U_m\|_2^2.
\]
Take $\varepsilon=\|\widetilde W_m-U_m\|_2^{2/3}$ when the norm
is nonzero; the zero case is immediate. The comparison error
decays polynomially and preserves the deterministic rate
$C/\ln(em)$. Averaging with \eqref{clock-logarithmic} proves
\eqref{kolmogorov-critical}.

Under only a second moment, indicator smoothing gives
\[
 d_K(U,G)\le C\inf_{\varepsilon>0}
   \{\varepsilon+\varepsilon^{-2}\zeta_2(U,G)\}
   \le C\zeta_2(U,G)^{1/3}.
\]
Theorem~\ref{thm:low-second-critical}, proved below independently
of the Kolmogorov bounds, therefore yields
\eqref{kolmogorov-second-critical}.
\end{manuscriptproof}

\begin{manuscriptproof}[Proof of Theorem~\ref{thm:zeta-low}]
In the diffusive regime, put $t=D_m/B_m\ge c$, so that
$U_m=t\widetilde W_m$ and $t^2-1=\Delta_m/B_m^2$.
For $f\in\Lambda_s$, subtract its affine part when $s>1$. Then
\[
 |f(tx)-f(x)|\le
 \begin{cases}
 |t-1|^s|x|^s,&0<s\le1,\\
 C_s|t^2-1||x|^s,&1<s\le2.
 \end{cases}
\]
Conditioning on the amplitudes and using
\eqref{conditional-self-moments} and \eqref{square-moment} gives
\begin{equation}
 \zeta_s(\widetilde W_m,U_m)
 \le C_s m^{-(s\wedge1)\rho/(2+\rho)},\qquad 0<s\le2.
\end{equation}
Combine this with \eqref{deterministic-low-zeta}.
At criticality, \eqref{critical-normalization-l2} and the uniform
second moments give
\[
 \zeta_s(\widetilde W_m,U_m)
 \le C_s\|\widetilde W_m-U_m\|_2^{s\wedge1}.
\]
The comparison error is polynomially small, so the deterministic
bound retains its logarithmic order.
Averaging over \eqref{mixture} and using
\eqref{clock-all-negative}--\eqref{clock-logarithmic} proves
\eqref{low-moment-zeta} and \eqref{zeta-critical}.
When $3-4p_r\ge1$, the contribution from $V_m^{-1}$ is absorbed
by the leading diffusive term.
Finally, \eqref{conditional-self-moments},
Theorem~\ref{thm:kolmogorov} and the weighted distribution-function
estimate \eqref{weighted-cdf-comparison} below, with $j=2$, give
\eqref{zeta-kolmogorov-consequence}.
\end{manuscriptproof}

\begin{manuscriptproof}[Proof of Theorem~\ref{thm:low-second-critical}]
Lemma~\ref{lem:critical-comparison}, \eqref{unit-low-zeta} and
\eqref{clock-logarithmic} give
\[
 \zeta_2(\widetilde W_m,G)
 \le C\{m/Q_m+m^{-1}+V_m^{-1}\}\le C/\ln(em),\qquad
 \zeta_2(W_n,G)\le C/\ln(en).
\]
Gaussian convolution yields
\[
 \zeta_s(U,G)\le C_s\inf_{\varepsilon>0}
   \{\varepsilon^s+\varepsilon^{s-2}\zeta_2(U,G)\}
   \le C_s\zeta_2(U,G)^{s/2},\qquad 0<s<2,
\]
which proves \eqref{low-second-zeta}.
\end{manuscriptproof}

\begin{manuscriptproof}[Proof of Theorem~\ref{thm:wasserstein-low}]
In the diffusive regime, put $t=D_m/B_m\ge c$. Conditioning on
the amplitudes as in \eqref{conditional-self-moments} gives,
for $1\le s\le2$,
\[
 \|\widetilde W_m-U_m\|_s^s\le C_s\E|t-1|^s
 \le C_s
 \begin{cases}
 m^{-s\rho/(2+\rho)},&s\le1+\rho/2,\\
 m^{-\rho/2},&s>1+\rho/2.
 \end{cases}
\]
The second case uses
$|\sqrt{1+u}-1|^s\le C_s|u|^{1+\rho/2}$ for $u\ge c-1$.
Consequently,
\begin{equation}
 \|\widetilde W_m-U_m\|_s
 \le C_s m^{-\min\{\rho/(2+\rho),\rho/(2s)\}}.
\end{equation}
Combine this with \eqref{deterministic-low-wasserstein}.
At criticality, use
\[
 \mathcal W_s(\widetilde W_m,U_m)
 \le\|\widetilde W_m-U_m\|_2,\qquad 1\le s\le2,
\]
and \eqref{critical-normalization-l2}; the comparison error
is polynomially small.

Now apply \eqref{transport-mixture} with
\eqref{clock-all-negative}--\eqref{clock-logarithmic}.
The leading diffusive clock exponents are below one except
$s\min\{\rho/(2+\rho),\rho/(2s)\}=1$ at $(\rho,s)=(2,2)$,
which gives the stated $\sqrt{\ln n}$ factor.
The $V_m^{-1}$ contribution is absorbed by the first term whenever
$3-4p_r\ge1$. This proves \eqref{low-moment-wasserstein}
and \eqref{wasserstein-critical}.
Finally, \eqref{conditional-self-moments},
Theorem~\ref{thm:kolmogorov} and
\eqref{weighted-cdf-comparison}, with $j=1$, give
\eqref{wasserstein-kolmogorov-consequence}.
\end{manuscriptproof}

\begin{manuscriptproof}[Proof of Corollary~\ref{cor:wasserstein-second-critical}]
For $1\le s\le2$, \cite[Theorem~3.1]{Rio2009} gives
\[
 \mathcal W_s(W_n,G)\le\mathcal W_2(W_n,G)
 \le C\zeta_2(W_n,G)^{1/2}.
\]
The assertion follows from Theorem~\ref{thm:low-second-critical}.
\end{manuscriptproof}

\subsection{Sharpness and distribution-function comparisons}

We use the following elementary consequence of subgaussian tails.
If $U,V$ satisfy
$\Pp(|U|>u)+\Pp(|V|>u)\le K e^{-cu^2}$ for all $u\ge0$,
with fixed $K,c>0$, and $\varepsilon=d_K(U,V)\in(0,1/2]$, then
\begin{equation}\label{weighted-cdf-comparison}
 \int_{\R}|x|^{j-1}|F_U(x)-F_V(x)|\,dx
 \le C_j\varepsilon\{\ln(e/\varepsilon)\}^{j/2},
 \qquad j=1,2,3.
\end{equation}
Indeed, split the integral at
$R=\{c^{-1}\ln(K'/\varepsilon)\}^{1/2}$ with a sufficiently
large fixed $K'$. The part on $[-R,R]$ is at most
$2\varepsilon R^j/j$, and the remaining part is bounded by
$2K\int_R^\infty x^{j-1}e^{-cx^2}\,dx$.
The case $\varepsilon=0$ follows by equality in distribution.
For $\varepsilon>1/2$, the analogous bound follows after changing
the constant. Integration by parts now gives the
$\mathcal W_1$ bound at $j=1$, the $\zeta_2$ bound for centered
variables at $j=2$, and the anchored third-order bound at $j=3$.
For the ideal metrics one subtracts the affine or quadratic
Taylor polynomial as appropriate; the remaining derivative is
bounded by $|x|$ or $x^2/2$, respectively. The full third-order
metric additionally requires equality of second moments.

\begin{manuscriptproof}[Proof of Proposition~\ref{prop:kolmogorov-sharp}]
If $p_r=1/2$, the unit-step embedded walk is a simple symmetric
random walk, so $\Pp(T_m=0)\asymp m^{-1/2}$ for even $m$.
Equations \eqref{clock-parity}, \eqref{negative-moment} and
\eqref{clock-positive}, with Markov's inequality for the upper
clock tail, show that some fixed $0<c_1<C_1<\infty$ satisfy
\[
 \Pp(c_1n^b\le\Sigma_n\le C_1n^b,\
                         \Sigma_n\ \hbox{is even})\ge c_2>0
\]
for all large $n$. Consequently,
\[
 \Pp(W_n=0)
 =\E\{\Pp(T_m=0)|_{m=\Sigma_n}\}
 \ge c n^{-b/2}.
\]
A jump of size $q$ in a distribution function forces its
Kolmogorov distance from any continuous distribution to be at
least $q/2$. This proves the lower bound in
\eqref{balanced-optimal-rate}; the upper bound follows from
Theorem~\ref{thm:kolmogorov}.
\end{manuscriptproof}

\begin{manuscriptproof}[Proof of Proposition~\ref{prop:normal-lower-bounds}]
Suppose first that $1/2<p_r<3/4$, and put
$a=2p_r-1\in(0,1/2)$ and $\alpha=1-2a\in(0,1)$.
By \eqref{unit-subgaussian},
\[
 0<\sum_{j=2}^\infty
       \frac{A_j^2\E T_{j-1}^2}{(j-1)^2}
 \le C\sum_{j=2}^\infty j^{-1-2a}<\infty.
\]
Thus \eqref{unit-variance-deficit} and $V_m\asymp m^\alpha$
give $1-\E T_m^2/Q_m\asymp m^{-\alpha}$ for $m\ge2$.
Averaging with \eqref{exact-negative} and
\eqref{clock-singleton}, and testing with $x^2/2\in\Lambda_2$, yields
\[
 \zeta_2(W_n,G)\ge\tfrac12(1-\E W_n^2)\ge c n^{-b\alpha}.
\]
For the matching upper bound, the unit-step estimate
\eqref{unit-low-zeta}, which is
\cite[Theorem~3.4]{Dedecker2023}, gives
\[
 \zeta_2(W_n,G)
 \le C\E\{\Sigma_n^{-1}+V_{\Sigma_n}^{-1}\}
 \le C\{n^{-b}\ln(en)+n^{-b\alpha}\}
 \le Cn^{-b\alpha}.
\]
At criticality, Theorem~\ref{thm:low-second-critical} gives the
upper bound. For the lower bound, \eqref{unit-variance-deficit} gives
\[
 0<\sum_{j=2}^\infty\frac{A_j^2\E T_{j-1}^2}{(j-1)^2}
 \le C\sum_{j=2}^\infty\frac{\ln(ej)}{j^2}<\infty,
 \qquad
 1-\E T_m^2/Q_m\asymp\{\ln(em)\}^{-1}\quad(m\ge2).
\]
Since $\Sigma_n\le n$ and $\Pp(\Sigma_n\ge2)\to1$, testing with
$x^2/2$ yields
\[
 \zeta_2(W_n,G)\ge\tfrac12|\E W_n^2-1|
 \ge\frac{c\Pp(\Sigma_n\ge2)}{\ln(en)}\ge\frac c{\ln n}.
\]
\end{manuscriptproof}

\subsection{Third-order bounds}

In the following estimates the test function $f$ satisfies
\begin{equation}\label{anchored-tests}
 |f''(x)-f''(y)|\le|x-y|,
 \qquad f(0)=f'(0)=f''(0)=0.
\end{equation}
In particular,
\begin{equation}\label{test-growth}
 |f''(x)|\le|x|,\quad |f'(x)|\le x^2/2,\quad
 |f(x)|\le|x|^3/6,\quad
 |f(tx)-f(x)|\le\tfrac16|t^3-1||x|^3\quad(t\ge0).
\end{equation}
The restriction in \eqref{anchored-tests} is used for auxiliary
estimates only. The full test class in \eqref{probability-distances} is
recovered after matching the first two moments.

We first record two truncation comparisons that will exploit
\eqref{conditional-self-moments}. Suppose that $U$ and $V$ satisfy
$\Pp(|U|>u)+\Pp(|V|>u)\le Ke^{-cu^2}$ for $u\ge0$, with fixed
$K,c>0$. Under any coupling, Cauchy--Schwarz and these tail bounds give
\[
 \E\{(|U|^j+|V|^j)\mathbf1_{\{|U|\vee|V|>R\}}\}
 \le C_j e^{-c'R^2},\qquad R\ge1,\quad j=2,3.
\]
Take an optimal $\mathcal W_1$ coupling and put
$\eta=\mathcal W_1(U,V)$. On $\{|U|\vee|V|\le R\}$,
\eqref{test-growth} bounds $|f(U)-f(V)|$ by
$CR^2|U-V|$, and $|U^2-V^2|\le2R|U-V|$.
Choosing $R$ to be a sufficiently large constant times
$\sqrt{\ln(e/\eta)}$ yields, for $0<\eta\le1/2$,
\begin{align}
 |\E U^2-\E V^2|&\le C\eta\sqrt{\ln(e/\eta)},
 \label{subgaussian-variance-comparison}\\
 \sup_{f\text{ satisfying }\eqref{anchored-tests}}
 |\E f(U)-\E f(V)|&\le C\eta\ln(e/\eta).
 \label{subgaussian-smooth-comparison}
\end{align}
Likewise, under an optimal $\mathcal W_2$ coupling,
$|U-V|^3\le2R|U-V|^2$ on that event. With
$\eta=\mathcal W_2(U,V)$, the same choice of $R$ gives
\begin{equation}\label{subgaussian-transport-comparison}
 \mathcal W_3(U,V)
 \le C\eta^{2/3}\{\ln(e/\eta)\}^{1/6},\qquad 0<\eta\le1/2.
\end{equation}
All three right-hand sides are interpreted as zero when $\eta=0$.
They are increasing functions of $\eta$ near zero, so an upper bound
for $\eta$ may be substituted in these estimates.

\begin{lemma}\label{lem:third-smooth}
Uniformly over the functions in \eqref{anchored-tests}, if
$0<p_r<3/4$ and \eqref{higher-moment} holds with $1\le\rho\le2$, then
\begin{equation}\label{fixed-third-moment}
 |\E f(\widetilde W_m)-\E f(G)|
   +|\E\widetilde W_m^2-1|
 \le C\{m^{-\rho/(2+\rho)}+V_m^{-1}\}.
\end{equation}
If $p_r=3/4$ and
$\E|Z_1-1|^3<\infty$, then
\begin{equation}\label{fixed-third-critical}
 |\E f(\widetilde W_m)-\E f(G)|
   +|\E\widetilde W_m^2-1|
 \le\frac{C}{\ln(em)}.
\end{equation}
\end{lemma}

\begin{manuscriptproof}
We first establish
\begin{equation}\label{unit-third-smooth}
 \left|\E f(T_m/\sqrt{Q_m})-\E f(G)\right|
 +\left|\E T_m^2/Q_m-1\right|\le C(m^{-1}+V_m^{-1}).
\end{equation}
For $u=aT_{j-1}/(j-1)$,
\[
 \E(\varepsilon_j^2\mid\mathcal F_{j-1})=1-u^2,
 \qquad
 \E(\varepsilon_j^3\mid\mathcal F_{j-1})=-2u(1-u^2).
\]
Smooth $f$ by $\N(0,m^{-1})$ and telescope
$A_j\varepsilon_j/\sqrt{V_m}$ against independent Gaussians of
variance $A_j^2/V_m$. At the $j$th replacement, the third
derivative is bounded by one and the fourth by
\[
 C\{m^{-1}+V_m^{-1}\textstyle\sum_{k=j+1}^m A_k^2\}^{-1/2}.
\]
Taylor's formula through degree three gives
\begin{align}
 &\left|\E f(T_m/\sqrt{Q_m})-\E f(G)\right|\nonumber\\
 &\quad\le C\left\{\frac1m
 +\frac{a^2}{V_m}\sum_{j=2}^m
                   \frac{A_j^2Q_{j-1}}{(j-1)^2}
 +\frac{|a|}{V_m^{3/2}}\sum_{j=2}^m
                   \frac{A_j^3\sqrt{Q_{j-1}}}{j-1}\right.\nonumber\\
 &\hspace{45mm}\left.
 +\frac1{V_m^2}\sum_{j=1}^m
   \frac{A_j^4}
   {\{m^{-1}+V_m^{-1}\sum_{k=j+1}^m A_k^2\}^{1/2}}
 \right\}.
\end{align}
Here \eqref{test-growth} bounds the smoothing cost by $C/m$;
\eqref{unit-subgaussian} gives the quadratic estimate
\[
 \E\left[T_{j-1}^2
       \left(1+\frac{|A_{j-1}T_{j-1}|}{\sqrt{V_m}}\right)\right]
 \le C Q_{j-1}.
\]
For $a<1/2$, $A_j\asymp j^{-a}$, $V_m\asymp m^{1-2a}$ and
$Q_j\asymp j$, whence
\begin{gather*}
 \frac{a^2}{V_m}\sum_{j=2}^m\frac{A_j^2Q_{j-1}}{(j-1)^2}
 +\frac{|a|}{V_m^{3/2}}\sum_{j=2}^m
                    \frac{A_j^3\sqrt{Q_{j-1}}}{j-1}
 \le C(m^{-1}+V_m^{-1}),\\
 \frac{C}{V_m^2}\sum_{j\le m/2}A_j^4
 +\frac{C}{m^2}\sum_{j>m/2}
              \left(\frac m{m-j+1}\right)^{1/2}
 \le C(m^{-1}+V_m^{-1}).
\end{gather*}
The second line bounds the fourth-order sum. For $a=1/2$,
$V_m\asymp\ln(em)$ and $Q_j\asymp j\ln(ej)$ give quadratic
and cubic bounds $C/V_m$ and $C/V_m^{3/2}$. Splitting the
fourth-order sum at $\sqrt m$ gives
\[
 \frac C{V_m^2}\left(\sum_{j\le\sqrt m}j^{-2}
                  +\sqrt m\sum_{j>\sqrt m}j^{-2}\right)
 \le C/V_m^2.
\]
Together with \eqref{unit-variance-deficit}, these estimates prove
\eqref{unit-third-smooth}.

Under a finite third moment, telescope the two-stage martingale
\eqref{two-stage} with variance proxies
$A_j^2/(A_m^2B_m^2)$ and $\sigma^2/B_m^2$, whose sum is one.
Taylor's formula through degree two and the same quadratic
estimate give
\begin{align}
 &|\E f(U_m)-\E f(G)|+|\E U_m^2-1|\nonumber\\
 &\quad\le C\left\{
 \frac{\sum_{j=1}^m A_j^3}{V_m^{3/2}}+\frac m{B_m^3}
 +\frac{a^2}{V_m}\sum_{j=2}^m
                  \frac{A_j^2Q_{j-1}}{(j-1)^2}\right\}
 \le C(m^{-1/2}+V_m^{-1}).
\end{align}
For $p_r<3/4$, condition on the amplitudes in
\eqref{conditional-self-moments} and use $t=D_m/B_m$ in
\eqref{test-growth}:
\begin{align}\label{third-denominator-smooth}
 |\E f(\widetilde W_m)-\E f(U_m)|
 &\le C\E\left|\left(1+\frac{\Delta_m}{B_m^2}\right)^{3/2}-1\right|
 \nonumber\\
 &\le C\left\{\frac{\E|\Delta_m|}{m}
             +\frac{\E|\Delta_m|^{3/2}}{m^{3/2}}\right\},\\
 |\E\widetilde W_m^2-\E U_m^2|
 &\le C\E|\Delta_m|/m.\nonumber
\end{align}
For $1\le\rho\le2$, \eqref{square-moment} gives
\[
 \E|\Delta_m|\le Cm^{2/(2+\rho)},\qquad
 \E|\Delta_m|^{3/2}\le Cm^{3/(2+\rho)},
\]
which proves \eqref{fixed-third-moment}. At criticality,
Lemma~\ref{lem:critical-comparison} with $s=3,2$ and
\eqref{unit-third-smooth} give \eqref{fixed-third-critical}, since
$m/Q_m+V_m^{-1}\le C/\ln(em)$.
\end{manuscriptproof}

\begin{manuscriptproof}[Proof of Theorem~\ref{thm:zeta-three}]
For $1\le\rho\le2$ in the diffusive regime, and for the
critical third-moment assertion, average Lemma~\ref{lem:third-smooth}
using \eqref{mixture} and
\eqref{clock-all-negative}--\eqref{clock-logarithmic}. The bounds
in \eqref{zeta-three-moment} and \eqref{zeta-three-critical} follow
for $|\E W_n^2-1|$ and the anchored expectation difference.
For $3-4p_r\ge1$, the $V_m^{-1}$ contribution is absorbed by
the leading term.

For $0<\rho<1$ and $p_r<3/4$, Theorem~\ref{thm:wasserstein-low} gives
$\mathcal W_1(W_n,G)\le C R_n(\rho)$.
The variables $W_n$ have uniform subgaussian tails by
\eqref{conditional-self-moments}. Apply
\eqref{subgaussian-variance-comparison}--\eqref{subgaussian-smooth-comparison}
with $U=W_n$ and $V=G$. Since $R_n(\rho)$ decays polynomially,
this bounds both the variance error and the anchored expectation
difference by $C R_n(\rho)\ln n$, as required for
\eqref{zeta-three-low-moment}.

For the additional estimate \eqref{zeta-three-kolmogorov},
Theorem~\ref{thm:kolmogorov} gives
$d_K(W_n,G)\le C\mathcal K_n(\rho)$. Uniform subgaussian tails
and \eqref{weighted-cdf-comparison}, with $j=2,3$, therefore give
\[
 |\E W_n^2-1|\le C\mathcal K_n(\rho)\ln n,\qquad
 \sup_{f\ {\rm satisfying}\ \eqref{anchored-tests}}
       |\E f(W_n)-\E f(G)|
 \le C\mathcal K_n(\rho)(\ln n)^{3/2}.
\]
These estimates also apply when $0<\rho<1$.

Equations \eqref{conditional-self-moments}, or
\eqref{unit-subgaussian}--\eqref{amplitude-noise-moments} at
criticality, imply $\sup_n\E|W_n|^3<\infty$. Hence
\[
 |\E f(\overline W_n)-\E f(W_n)|
 \le\tfrac16\E|W_n|^3|(\E W_n^2)^{-3/2}-1|
 \le C|\E W_n^2-1|.
\]
For $f\in\Lambda_3$, subtract its Taylor polynomial of degree
two. Its expectations agree under $\overline W_n$ and $G$,
since
\[
 \E\overline W_n=\E G=0,\qquad
 \E\overline W_n^2=\E G^2=1.
\]
Taking the supremum completes the proof. The variance correction
is applied after averaging over the clock.
\end{manuscriptproof}

\begin{manuscriptproof}[Proof of Theorem~\ref{thm:wasserstein-three}]
For martingale differences $\xi_j$ adapted to $(\mathcal H_j)$ with
$\sum_j\E\xi_j^2=1$, Guo's inequality
\cite[Theorem~8 in the preprint version]{Guo2026} gives
\begin{equation}\label{guo-third-bound}
 \mathcal W_3\left(\sum_j\xi_j,G\right)
 \le C\left\{
 (\sum_j\E|\xi_j|^3)^{1/3}
 +\left\|\sum_j\E(\xi_j^2\mid\mathcal H_{j-1})-1
       \right\|_{3/2}^{1/2}\right\}.
\end{equation}
For \eqref{two-stage}, \eqref{two-stage-bracket} and
\eqref{unit-subgaussian} bound both the $L^{3/2}$ norm and the
expectation of the variance deficit by
\[
 \frac{Ca^2}{V_m}\sum_{j=2}^m
             \frac{A_j^2Q_{j-1}}{(j-1)^2}.
\]
Normalize by the actual variance in \eqref{guo-third-bound}
and restore the proxy $B_m^2$ to obtain
\begin{align}\label{deterministic-wasserstein-three}
 \mathcal W_3(U_m,G)
 &\le C\left\{
 \frac{(\sum_{j=1}^m A_j^3)^{1/3}}{\sqrt{V_m}}
 +\frac{m^{1/3}}{B_m}\right.\nonumber\\
 &\hspace{12mm}\left.
 +\left(\frac{a^2}{V_m}\sum_{j=2}^m
                  \frac{A_j^2Q_{j-1}}{(j-1)^2}\right)^{1/2}
 \right\}.
\end{align}
The variance ratios tend to one. For $p_r<3/4$, conditioning
on the amplitudes as in \eqref{third-denominator-smooth} gives
\begin{align}
 \|\widetilde W_m-U_m\|_3^3
 &\le C\E|D_m/B_m-1|^3\nonumber\\
 &\le C\frac{\E|\Delta_m|^{3/2}}{B_m^3}
 \le C m^{-1/2}.
\end{align}
We used $|\sqrt{1+u}-1|^3\le|u|^{3/2}$ and
\eqref{square-moment} with $\rho=1$. Also,
\[
 \sum_{j=1}^m A_j^3\asymp
 \begin{cases}
 m^{4-6p_r},&p_r<2/3,\\
 \ln(em),&p_r=2/3,\\
 1,&p_r>2/3.
 \end{cases}
\]
The quadratic-variation term in
\eqref{deterministic-wasserstein-three} is bounded by
$C(m^{-1}+V_m^{-1})^{1/2}$. Consequently,
\begin{equation}
 \mathcal W_3(\widetilde W_m,G)\le C
 \begin{cases}
 m^{-1/6},&0<p_r<2/3,\\
 m^{-1/6}\{\ln(em)\}^{1/3},&p_r=2/3,\\
 m^{-(3-4p_r)/2},&2/3<p_r<3/4.
 \end{cases}
\end{equation}

At criticality, apply \eqref{guo-third-bound} to the directional
martingale to obtain
\[
 \mathcal W_3(T_m/\sqrt{Q_m},G)\le C\{\ln(em)\}^{-1/2}.
\]
To compare its statistic with $\widetilde W_m$, observe that
\begin{align*}
 \|\widetilde W_m-T_m/\sqrt{Q_m}\|_3
 &\le \frac{\|\widetilde S_m-T_m\|_3}{\sqrt{Q_m}}\\
 &\quad+
 \|T_m/\sqrt{Q_m}\|_3
 \left\|1-\left(1+\frac{D_m^2-Q_m}{Q_m}\right)^{-1/2}\right\|_3\\
 &\le C\left(\frac{m}{Q_m}\right)^{1/2}
 \le C\{\ln(em)\}^{-1/2}.
\end{align*}
Here $T_m\perp(D_m^2-Q_m)$,
$|1-(1+u)^{-1/2}|\le\sqrt u$ and
$\E(D_m^2-Q_m)^{3/2}\le Cm^{3/2}$.
Now \eqref{transport-mixture} and
\eqref{clock-all-negative}--\eqref{clock-logarithmic} prove
\eqref{wasserstein-three-moment}; the diffusive clock exponents
are $1/2$, $1/2$ and $3(3-4p_r)/2$, all below one.

Finally, if $\E|Z_1-1|^4<\infty$ and $0<p_r\le1/2$,
Theorem~\ref{thm:wasserstein-low} gives
$\mathcal W_2(W_n,G)\le Cn^{-b/2}\sqrt{\ln n}$, since
$3-4p_r\ge1$. The uniform subgaussian bound
\eqref{conditional-self-moments} permits the use of
\eqref{subgaussian-transport-comparison}. Substitution gives
\[
 \mathcal W_3(W_n,G)
 \le C\{n^{-b/2}\sqrt{\ln n}\}^{2/3}(\ln n)^{1/6}
 =Cn^{-b/3}\sqrt{\ln n},
\]
which proves \eqref{wasserstein-three-fourth}.

\end{manuscriptproof}

\begin{manuscriptproof}[Proof of Remark~\ref{rem:intermediate-orders}]
Let $2<s<3$.
For $f\in\Lambda_s$, let $f_\varepsilon(x)=\E f(x+\varepsilon G)$.
The H\"older condition on $f''$ and Gaussian integration by parts give
\[
 \|(f-f_\varepsilon)''\|_\infty\le C_s\varepsilon^{s-2},\qquad
 \|f_\varepsilon'''\|_\infty\le C_s\varepsilon^{s-3}.
\]
Thus, for variables $U,V$ with equal first two moments,
\[
 \zeta_s(U,V)\le C_s\{
 \varepsilon^{s-2}\zeta_2(U,V)+
 \varepsilon^{s-3}\zeta_3(U,V)\}.
\]
Optimizing in $\varepsilon>0$ yields
\begin{equation}\label{zeta-interpolation}
 \zeta_s(U,V)\le C_s
 \zeta_2(U,V)^{3-s}\zeta_3(U,V)^{s-2}.
\end{equation}
The zero-distance cases follow by taking limits.
After subtracting the affine part of a test in $\Lambda_2$, a
scaling comparison gives
\[
 \zeta_2(\overline W_n,G)
 \le\zeta_2(W_n,G)+C|\E W_n^2-1|
\]
for all sufficiently large $n$.
Theorems~\ref{thm:zeta-low} and~\ref{thm:zeta-three}, followed
by \eqref{zeta-interpolation}, prove \eqref{intermediate-zeta}.
\end{manuscriptproof}

\begin{manuscriptproof}[Proof of Remark~\ref{rem:intermediate-wasserstein}]
For \eqref{intermediate-wasserstein}, the quantile coupling on
$\R$ is optimal simultaneously at orders $2$, $s$ and $3$.
Apply H\"older's inequality with
$1/s=\theta/2+(1-\theta)/3$, where $\theta=2(3-s)/s$.
With a finite third moment and $0<p_r\le5/8$, the endpoint bounds
are $\mathcal W_2(W_n,G)\le Cn^{-b/4}$ and
$\mathcal W_3(W_n,G)\le Cn^{-b/6}$, which give the stated example.
\end{manuscriptproof}

\section{Proofs of the moderate-deviation results}\label{sec:md-proofs}

\begin{manuscriptproof}[Proof of Theorem~\ref{thm:expsquare}]
Let $p_r<3/4$, fix $c>0$, and restrict to
$0\le x\le c\sqrt{\ln n}$. Choose $L>3/(2\kappa)$ and put
$M_n=\lceil(\ln n)^L\rceil$. The right-hand side of
\eqref{fixed-diff} is $o(1)$ uniformly for $M_n\le m\le n$,
since it is bounded by $C(\ln n)^{3/2-L\kappa}$.
The range condition $x=o(m^\kappa)$ is also satisfied uniformly.
Exponentiating \eqref{fixed-diff} therefore gives
\[
 \left|\frac{\Pp(\widetilde W_m\ge x)}{\Phib(x)}-1\right|
 \le C\{x^3+(1+x)\ln m\}m^{-\kappa}.
\]
The large-clock contribution is at most
\begin{equation}
 C\{x^3+(1+x)\ln n\}\E\Sigma_n^{-\kappa}
 \le C\{x^3+(1+x)\ln n\}n^{-b\kappa}.
\end{equation}
For the complement, \eqref{negative-moment} gives
$\Pp(\Sigma_n<M_n)\le Cn^{-b}(\ln n)^L$.
Mills' inequality,
\begin{equation}\label{mills}
 \Phib(x)^{-1}\le C(1+x)e^{x^2/2},\qquad x\ge0,
\end{equation}
therefore yields, after averaging over \eqref{mixture},
\begin{align}\label{diffusive-tail-two-errors}
 \left|\frac{\Pp(W_n\ge x)}{\Phib(x)}-1\right|
 &\le C\{x^3+(1+x)\ln n\}n^{-b\kappa}\nonumber\\
 &\quad+C(1+x)e^{x^2/2}n^{-b}(\ln n)^L.
\end{align}
If $c<\sqrt{2b}$, the second term is at most
$Cn^{-b+c^2/2}(\ln n)^{L+1/2}=o(1)$; the first term also tends
to zero uniformly. This proves \eqref{bounded-diff-range}.
If $c<\sqrt{2b(1-\kappa)}$, the second term is
$o(n^{-b\kappa})$ and is absorbed by the first.
Taking the logarithm of the ratio then proves \eqref{bounded-diff}.
Reflection gives both assertions for the lower tail.

For $p_r=3/4$, \eqref{negative-moment} gives
\begin{equation}
 \Pp(\Sigma_n<n^{b/2})\le Cn^{-b/2}.
\end{equation}
Put $H_n(x)=\{x^3+(1+x)\ln\ln n\}/\sqrt{\ln n}$.
For $n^{b/2}\le m\le n$, \eqref{fixed-critical} gives
\[
 e^{-CH_n(x)}\le
 \frac{\Pp(\widetilde W_m\ge x)}{\Phib(x)}\le e^{CH_n(x)}.
\]
Consequently, uniformly for $x=o(\sqrt{\ln n})$,
\[
 (1-Cn^{-b/2})e^{-CH_n(x)}
 \le\frac{\Pp(W_n\ge x)}{\Phib(x)}
 \le e^{CH_n(x)}+n^{-b/2+o(1)}.
\]
Since $H_n(x)\ge\ln\ln n/\sqrt{\ln n}$, this proves
\eqref{bounded-critical}. If $x=o((\ln n)^{1/6})$, then
$H_n(x)=o(1)$.
\end{manuscriptproof}

\begin{manuscriptproof}[Proof of Corollary~\ref{cor:moment}]
By continuity of $\Phi$ and \eqref{mills},
\[
 \left|\frac{\Pp(W_n\ge x)}{\Phib(x)}-1\right|
 \le C(1+x)e^{x^2/2}d_K(W_n,G).
\]
Theorem~\ref{thm:kolmogorov} gives
\eqref{moment-diff}, \eqref{moment-critical} and
\eqref{moment-second-critical}; \eqref{balanced-kolmogorov}
gives the simplification for $p_r\le5/8$. At the stated endpoints,
\[
 \begin{aligned}
 &(1+\sqrt{2\gamma\ln n})n^\gamma
   \mathcal K_n(\rho)\longrightarrow0,\\
 &\{1+\sqrt{2(1-\varepsilon)\ln\ln n}\}
    (\ln n)^{-\varepsilon}\longrightarrow0,\\
  &(1+\sqrt{2\gamma\ln\ln n})(\ln n)^{\gamma-1/3}
    \longrightarrow0\quad\text{if }0<\gamma<1/3.
 \end{aligned}
\]
Reflection gives the lower tails.
\end{manuscriptproof}

\section{Proofs of the large-deviation results}\label{sec:large-proofs}

\begin{lemma}\label{lem:embedded-pld}
Suppose that $0<p_r<1$. Under \eqref{regular-tail}, if
$h_m/\sqrt{Q_m\ln(em)}\to\infty$, then
\begin{equation}
 \sup_{x\ge h_m}
 \left|\frac{\Pp(\pm\widetilde S_m>x)}
 {m\overline F(x)/2}-1\right|\longrightarrow0.
\end{equation}
\end{lemma}

\begin{manuscriptproof}
Put $R_k=\sum_{j=1}^k(\widetilde Z_j-1)$, $R_0=0$, and let
$(R_k')$ be an independent copy. The independent-sum estimate in
\cite[Theorem~8.1]{DenisovDiekerShneer2008}, including its
fixed-$k$ counterpart, gives
\begin{equation}\label{iid-pld-uniform}
 \sup_{1\le k\le m}\sup_{x\ge h_m}
 \left|\frac{\Pp(R_k>x)}{k\overline F(x)}-1\right|\to0,
 \qquad h_m\gg\sqrt{m\ln(em)}.
\end{equation}
Choose $\varepsilon_m\downarrow0$ such that
\[
 \varepsilon_m\ge m^{-1/4},\qquad
 \frac{\varepsilon_m^2h_m^2}{Q_m\ln(em)}\longrightarrow\infty.
\]
By \eqref{unit-subgaussian} and the one-sided Bernstein inequality
($1-Z_1\le1$),
\begin{align*}
 \Pp(|T_m|>v)&\le2e^{-v^2/(2Q_m)},\\
 \sup_{0\le k\le m}\Pp(-R_k>v)
 &\le\exp\{-cv^2/(m+v)\}.
\end{align*}
For $h_m\le x\le m^2$, regular variation therefore gives
\begin{equation}\label{pld-negligible-errors}
 \Pp(|T_m|>\varepsilon_mx)
 +\sup_{k\le m}\Pp(-R_k>\varepsilon_mx)
 =o(m\overline F(x)).
\end{equation}
The events in \eqref{pld-negligible-errors} are empty for $x>m^2$
and large $m$. Also,
\[
 \sup_{k\le m}\Pp(R_k>\varepsilon_mx)
 \le\frac{m\sigma^2}{\varepsilon_m^2x^2}=o(1).
\]
Conditional on $T_m=2k-m$,
\[
 \widetilde S_m\eqdist2k-m+R_k-R_{m-k}',
 \qquad R_k\ \hbox{and}\ R_{m-k}'\ \hbox{independent}.
\]
For $|2k-m|\le\varepsilon_mx$, the conditional tail is between
\[
 \Pp(R_k>(1+2\varepsilon_m)x)
 \{1-\Pp(R_{m-k}'>\varepsilon_mx)\}
\]
and
\[
 \Pp(R_k>(1-2\varepsilon_m)x)
       +\Pp(-R_{m-k}'>\varepsilon_mx).
\]
Equations \eqref{iid-pld-uniform}--\eqref{pld-negligible-errors},
regular variation and $\E T_m=0$ yield, uniformly for $x\ge h_m$,
\[
 \Pp(\widetilde S_m>x)
 =\overline F(x)\left\{
 \E\left(\frac{m+T_m}{2};\ |T_m|\le\varepsilon_mx\right)
 +o(m)\right\}
 \sim\frac m2\overline F(x).
\]
Reflection gives the other sign.
\end{manuscriptproof}

\begin{manuscriptproof}[Proof of Theorem~\ref{thm:position-pld}]
Choose $C$ sufficiently large and set
$M_n=\lfloor n^b(C\ln n)^{1-b}\rfloor$. By
\eqref{clock-upper-ld}--\eqref{clock-small-rate} and
\eqref{clock-positive},
\begin{equation}
 \Pp(\Sigma_n>M_n)\le n^{-2\gamma-3},\qquad
 \E(\Sigma_n;\Sigma_n>M_n)
 \le(\E\Sigma_n^2)^{1/2}\Pp(\Sigma_n>M_n)^{1/2}
 =o(n^b).
\end{equation}
Moreover,
\[
 \sup_{m\le M_n}\sqrt{Q_m\ln(em)}
 \le C
 \begin{cases}
 n^{b/2}(\ln n)^{1-b/2},&p_r<3/4,\\
 n^{b/2}(\ln n)^{(3-b)/2},&p_r=3/4,\\
 n^{b(2p_r-1)}(\ln n)^{(2p_r-1)(1-b)+1/2},&3/4<p_r<1.
 \end{cases}
\]
Lemma~\ref{lem:embedded-pld} and the finite-convolution asymptotic
for bounded $m$ imply
\[
 \sup_{1\le m\le M_n}\sup_{x\ge h_n}
 \left|\frac{\Pp(\widetilde S_m>x)}{m\overline F(x)/2}-1\right|
 \longrightarrow0.
\]
For $h_n\le x\le n^2$, \eqref{position-mixture} now gives
\begin{align*}
 \Pp(S_n>x)
 &=\frac{\overline F(x)}2
       \E(\Sigma_n;\Sigma_n\le M_n)(1+o(1))
       +O(n^{-2\gamma-3}),\\
 \E\Sigma_n\,\overline F(x)
 &\ge c n^b\overline F(n^2)
 \ge c n^{b-2\gamma-1}.
\end{align*}
Thus \eqref{position-pld} holds uniformly on this interval.
For $x\ge n^2$, Lemma~\ref{lem:embedded-pld} applies uniformly
over every $1\le m\le n$; averaging gives the same result.
The third line of the normalizer bound follows from
\eqref{superdiffusive-normalizer}; the embedded argument uses only
\eqref{unit-subgaussian}, not a diffusive or critical central limit
theorem. Finally, $h_n^*=o(n)$ throughout $0<p_r<1$.
Thus \eqref{clock-positive} and $\ell(nx)/\ell(n)\to1$ prove
\eqref{position-ballistic-pld} throughout this range.
\end{manuscriptproof}

\begin{manuscriptproof}[Proof of Proposition~\ref{prop:heavy-threshold}]
Put $\nu=0$ if $p_r<3/4$ and $\nu=1$ if $p_r=3/4$.
Then $Q_m\asymp m\{\ln(em)\}^{\nu}$ and
\[
 h_n^*=n^{b/2}(\ln n)^{(2-b+\nu)/2}.
\]
Given $h_n=o(h_n^*)$, set
\begin{equation}\label{threshold-test-sequence}
 x_n=\max\left\{h_n,\frac{h_n^*}{(\ln n)^{(2-b)/4}}\right\},
 \qquad
 k_n=\left\lceil
       \left\{\frac{x_n^2}{n^b(\ln n)^\nu}\right\}^{1/(2-b)}
      \right\rceil.
\end{equation}
These choices satisfy
\[
 x_n\ge h_n,\qquad x_n=o(h_n^*),\qquad
 (\ln n)^{1/2}\le k_n=o(\ln n),\qquad
 \ln x_n=\tfrac b2\ln n+O(\ln\ln n).
\]
Let
\[
 \mathcal B_n=\{m\in\mathbb N:
        c_0n^bk_n^{1-b}\le m\le C_0n^bk_n^{1-b}\},
\]
with the constants of Lemma~\ref{lem:clock-band}. That lemma gives
$\Pp(\Sigma_n\in\mathcal B_n)\ge e^{-Ck_n}$.
Uniformly for $m\in\mathcal B_n$, we have
\begin{equation}\label{threshold-band-scales}
 \ln m\sim b\ln n,\qquad
 \frac{x_n^2}{Q_m}\asymp k_n,\qquad
 \frac{m}{x_n^2}\le\frac{C}{k_n(\ln n)^\nu}=o(1).
\end{equation}

Apply Lemma~\ref{lem:fixed-cmd} to the unit-step walk, for which
$\widetilde W_m=T_m/\sqrt{Q_m}$. Its argument
$2x_n/\sqrt{Q_m}$ is of order $\sqrt{k_n}$.
In the diffusive regime, the error in \eqref{fixed-diff} is
$o(k_n)$, since $m$ grows polynomially in $n$.
At criticality, $\sqrt{k_n}=o(\sqrt{\ln m})$, and the error in
\eqref{fixed-critical} is bounded by
\[
 C\frac{k_n^{3/2}+(1+\sqrt{k_n})\ln\ln n}{\sqrt{\ln n}}
 =o(k_n).
\]
Mills' inequality consequently yields, uniformly on this clock band,
\begin{equation}\label{direction-band-lower}
 \Pp(T_m>2x_n)\ge e^{-Ck_n}.
\end{equation}
To obtain the strict inequality, one may use the next allowed
lattice point above $2x_n$; its distance from $2x_n$ is at most
two and does not change these estimates.

Write $N_m=\widetilde S_m-T_m=\sum_jX_j(\widetilde Z_j-1)$.
Conditional on the complete direction vector,
\[
 \E(N_m\mid X_1,\ldots,X_m)=0,\qquad
 \E(N_m^2\mid X_1,\ldots,X_m)=m\sigma^2.
\]
Conditional Chebyshev's inequality and
\eqref{threshold-band-scales} imply
\begin{align*}
 \Pp(\widetilde S_m>x_n)
 &\ge \E\left[
       \mathbf1_{\{T_m>2x_n\}}
       \Pp(N_m\ge-x_n\mid X_1,\ldots,X_m)\right]\\
 &\ge(1-m\sigma^2/x_n^2)\Pp(T_m>2x_n)
 \ge\tfrac12e^{-Ck_n}.
\end{align*}
The use of the conditional bound is important: it multiplies the
small directional probability rather than subtracting an
unconditional polynomial error from it.
Averaging on $\mathcal B_n$ gives
\begin{equation}\label{position-band-lower}
 \Pp(S_n>x_n)\ge e^{-Ck_n}=n^{-o(1)}.
\end{equation}
On the other hand, regular variation,
$\ln x_n=(b/2)\ln n+O(\ln\ln n)$ and \eqref{clock-positive} give
\[
 \tfrac12\E\Sigma_n\,\overline F(x_n)
 =n^{-b(\gamma-2)/2+o(1)}.
\]
Since $\gamma>2$, this and \eqref{position-band-lower} prove
\eqref{heavy-threshold-necessity} for the positive tail.
Reflection of the direction sequence gives the negative tail.
\end{manuscriptproof}

\begin{manuscriptproof}[Proof of Theorem~\ref{thm:position-weibull}]
Fix $x>0$. Independence and \eqref{clock-singleton} give
\begin{equation}
 \Pp(S_n>nx)\ge\tfrac12\Pp(\Sigma_n=1)\overline F(nx),
 \qquad
 \liminf_{n\to\infty}n^{-\alpha}\ln\Pp(S_n>nx)
 \ge-\beta x^\alpha.
\end{equation}
For $0<\delta<\beta$, put $t_n=(\beta-\delta)(nx)^{\alpha-1}$.
Condition \eqref{weibull-ln-tail} implies
\[
 \E\{Z_1^2e^{(\beta-\delta)Z_1^\alpha}\}<\infty,
 \qquad t_nz\le(\beta-\delta)z^\alpha\quad(0<z\le nx).
\]
Taylor's formula and $\E Z_1=1$ yield the truncated expectation bound
\begin{equation}\label{weibull-truncated-mgf}
 \E(e^{t_nZ_1};Z_1\le nx)
 \le1+t_n+\frac{t_n^2}{2}
                  \E\{Z_1^2e^{(\beta-\delta)Z_1^\alpha}\}
 \le\exp\{t_n+Ct_n^2\}.
\end{equation}
Let $0<\varepsilon<\min\{x,1\}$ and $k_n=\lceil\varepsilon n\rceil$.
Since the amplitudes are positive, \eqref{position-mixture},
\eqref{clock-upper-ld} and \eqref{weibull-truncated-mgf} give
\begin{align*}
 \Pp(S_n>nx)
 &\le\Pp(\Sigma_n>\varepsilon n)
         +\Pp\left(\sum_{j=1}^{k_n}\widetilde Z_j>nx\right)\\
 &\le b^{-1}e^{-nJ_b(\varepsilon)}
       +k_n\overline F(nx)
       +\exp\{-t_nnx+k_n(t_n+Ct_n^2)\}.
\end{align*}
Here $nt_n^2=o(n^\alpha)$ and $J_b(\varepsilon)>0$. Hence
\[
 \limsup_{n\to\infty}n^{-\alpha}\ln\Pp(S_n>nx)
 \le-(\beta-\delta)x^{\alpha-1}(x-\varepsilon).
\]
Let $\varepsilon,\delta\downarrow0$ and use reflection to obtain
\eqref{position-weibull-tail}. Furthermore,
\[
 \frac{\E|S_n|}{n}\le\frac{\E\Sigma_n}{n}\to0,
 \qquad
 \lim_{R\to\infty}\limsup_{n\to\infty}n^{-\alpha}
       \ln\Pp(|S_n|>nR)=-\infty.
\]
For $0<\varepsilon<x$, subtraction of the tails gives
\[
 \liminf_{n\to\infty}n^{-\alpha}
 \ln\Pp(x-\varepsilon<S_n/n<x+\varepsilon)
 \ge-\beta(x-\varepsilon)^\alpha.
\]
Reflection gives the negative half-line, and $S_n/n\to0$ in
probability gives the lower bound at zero. For a closed set $B$,
\[
 \Pp(S_n/n\in B)\le\Pp(|S_n|/n\ge\inf_{y\in B}|y|)
\]
gives the upper bound, by \eqref{position-weibull-tail} and
continuity of $I$. The argument uses no restriction on $p_r$.
\end{manuscriptproof}

\begin{manuscriptproof}[Proof of Theorem~\ref{thm:position-light}]
The numbers of positive and negative embedded directions are
$(m+T_m)/2$ and $(m-T_m)/2$. Consequently,
\begin{equation}\label{marked-mgf-identity}
 \E e^{t\widetilde S_m}
 =\{M(t)M(-t)\}^{m/2}
 \E\exp\left\{\frac{T_m}{2}\ln\frac{M(t)}{M(-t)}\right\},
 \qquad
 \lim_{m\to\infty}m^{-1}\ln\E e^{t\widetilde S_m}=\psi(t).
\end{equation}
The urn function for $(m+T_m)/2$ is
$u\mapsto1-p_r+(2p_r-1)u$. For every $0<p_r<1$, both endpoint
values lie strictly between zero and one, and the unique fixed
point is $1/2$. Thus the linear urn results apply throughout this
parameter range, not only in the diffusive and critical regimes. By
\cite[Theorem~9 and Corollary~12]{Franchini2017}, with the
correction in \cite{Franchini2025},
$\lambda_{p_r}\in C^1(\R)$; for $p_r=1/2$,
$\lambda_{1/2}(t)=\ln\cosh t$. Thus
\[
 \psi\in C^1(\R),\quad\psi(t)\ge0,\quad
 L_b\in C^1(\R),\quad L_b'(0)=0,\quad
 \Lambda=L_b\circ\psi\in C^1(\R).
\]
For fixed $t$ and $\varepsilon>0$, \eqref{marked-mgf-identity}
provides $C_{t,\varepsilon}\ge1$ such that, for all $m\ge1$,
\[
 C_{t,\varepsilon}^{-1}e^{m(\psi(t)-\varepsilon)}
 \le\E e^{t\widetilde S_m}
 \le C_{t,\varepsilon}e^{m(\psi(t)+\varepsilon)}.
\]
Averaging over \eqref{position-mixture} and using
\eqref{clock-mgf-limit},
\begin{align*}
 L_b(\psi(t)-\varepsilon)
 &\le\liminf_{n\to\infty}n^{-1}\ln\E e^{tS_n}\\
 &\le\limsup_{n\to\infty}n^{-1}\ln\E e^{tS_n}
 \le L_b(\psi(t)+\varepsilon).
\end{align*}
Letting $\varepsilon\downarrow0$ proves
\eqref{position-cumulant-limit}. Since $\Lambda$ is finite and
differentiable on $\R$, the G\"artner--Ellis theorem
\cite[Theorem~2.3.6]{DemboZeitouni1998} gives
\eqref{position-light-rate} and \eqref{position-ldp}.
\end{manuscriptproof}

\begin{manuscriptproof}[Proof of Theorem~\ref{thm:position-sharp}]
Use $H_n(w)=\E w^{\Sigma_n}$ and
\eqref{clock-generating-function}.
The real-parameter leading term below is already contained in
\cite{BercuFavaro2026}; the pole argument also provides the local
complex-parameter extension needed for exponential tilting.
For real $w>1$, the unique pole in $|z|<1$ is
\[
 \rho(w)=1-(1-w^{-1})^{1/b}\in(0,1),
 \qquad
 \mathop{\rm Res}_{z=\rho(w)}
 \frac{(1-z)^{b-1}}{(1-z)^b-1+w^{-1}}=-\frac1b.
\]
Removing this pole and applying Cauchy's coefficient formula gives,
locally uniformly for complex $w$ near each real $w>1$,
\begin{equation}
 H_n(w)=b^{-1}\rho(w)^{-n}\{1+O(e^{-cn})\}.
\end{equation}

For balanced directions the embedded increments are independent
with the law of a fair sign times $Z_1$. Put
$g(z)=\E\cosh(zZ_1)$. For fixed $t>0$,
\begin{equation}\label{sharp-complex-mgf}
 \E e^{zS_n}=H_n(g(z))
 =b^{-1}e^{n\Lambda(z)}\{1+O(e^{-cn})\}
\end{equation}
locally uniformly near $z=t$, with the analytic continuation of
$\Lambda$ at $t$.
The density assumption and two integrations by parts imply
\begin{equation}\label{sharp-fourier-gap}
 \sup_{|u|\ge\delta}|g(t+iu)|<g(t)\quad(\delta>0),
 \qquad |g(t+iu)|\le C_t(1+|u|)^{-2}.
\end{equation}
Moreover, $\Lambda''(t)>0$, $\Lambda'(0)=0$ and
$\Lambda'(t)\to\sup\operatorname{supp}(Z_1)$ as $t\to\infty$,
so $t_x$ exists uniquely.

Define the tilted law by
\[
 \frac{d\Pp_t}{d\Pp}
 =\frac{e^{tS_n}}{\E e^{tS_n}},
 \qquad
 \chi_{n,t}(u)=\frac{H_n(g(t+iu))}{H_n(g(t))}.
\]
Equations \eqref{clock-mgf-upper}, \eqref{sharp-complex-mgf}
and \eqref{sharp-fourier-gap} give, for suitable
$0<\delta<R<\infty$,
\begin{equation}\label{sharp-fourier-bounds}
 |\chi_{n,t}(u)|\le
 \begin{cases}
 Ce^{-cnu^2},&|u|\le\delta,\\
 Ce^{-cn},&\delta\le|u|\le R,\\
 Ce^{-n\Lambda(t)}(1+|u|)^{-2},&|u|\ge R.
 \end{cases}
\end{equation}
For the last bound use $H_n(v)\le v$ for $0\le v\le1$.
Taking $t=t_x$, Taylor expansion in \eqref{sharp-complex-mgf}
and \eqref{sharp-fourier-bounds} yield
\[
 \int_{\R}\left|
 \E_t\exp\left\{iu\frac{S_n-nx}{\sqrt{n\Lambda''(t)}}\right\}
 -e^{-u^2/2}\right|\,du\longrightarrow0.
\]
If $f_{n,t}$ is the density of $S_n$ under $\Pp_t$, Fourier
inversion consequently gives
\begin{equation}
 \sup_{y\in\R}\left|
 \sqrt{n\Lambda''(t)}\,
 f_{n,t}\!\left(nx+y\sqrt{n\Lambda''(t)}\right)-\phi(y)
 \right|\longrightarrow0,
 \qquad \phi=\Phi'.
\end{equation}
Finally, by dominated convergence,
\begin{align*}
 \Pp(S_n\ge nx)
 &=e^{-tnx}\E e^{tS_n}
                  \int_0^\infty e^{-ty}f_{n,t}(nx+y)\,dy\\
 &\sim\frac{e^{-n(tx-\Lambda(t))}}{b\sqrt{2\pi n\Lambda''(t)}}
                  \int_0^\infty e^{-ty}\,dy
 =\frac{e^{-nI(x)}}{b\,t_x\sqrt{2\pi n\Lambda''(t_x)}}.
\end{align*}
To obtain the asserted uniformity, let $K$ be a compact subset of
$(0,\sup\operatorname{supp}(Z_1))$. The tilts $t_x$, $x\in K$,
range over a compact subset of $(0,\infty)$, on which
$\Lambda''$ is bounded above and bounded away from zero.
The pole expansion, Fourier bounds and Taylor expansion above hold
uniformly on this compact set of tilts. Consequently, the local
limit estimate is uniform in $x\in K$ as well as in its spatial
argument. After multiplication by $\sqrt{n\Lambda''(t_x)}$,
the last integrands are bounded by $Ce^{-t_*y}$, where
$t_* =\min_{x\in K}t_x>0$.
Uniform dominated convergence proves the compact-uniform
asymptotic equivalence. Reflection completes the proof.
\end{manuscriptproof}






\end{document}